\documentclass[a4paper,11pt, english]{smfart}
\usepackage{vmargin,graphicx,amsmath,amssymb,latexsym, enumerate, color}
\usepackage[french]{babel}
\newtheorem{thm}{Theorem}[section]  \newtheorem{cor}[thm]{Corollary}
\newtheorem{lem}[thm]{Lemma}   
  \newtheorem{prop}[thm]{Proposition}

\def\remark{\refstepcounter{thm}\bigskip\noindent\bf Remark \thethm\rm\ }

\newcommand{\preuve}[1][\!\!]{\bigskip\noindent{\it Proof #1. \ \ }}
\def\fin{\hfill$\Box$\\}
\newcounter{numexo}\newcounter{numq}
\newcounter{numsq}

\def\ccc{{\mathcal C}}
\def\eee{{\mathcal E}}\def\fff{{\mathcal F}}\def\ggg{{\mathcal G}} \def\hhh{{\mathcal H}}
 \def\kkk{{\mathcal K}}\def\lll{{\mathcal L}}
\def\mmm{{\mathcal M}} 
\def\sss{{\mathcal S}}

\def\norm#1{\left\Vert#1\right\Vert}

\def\abs#1{\left\vert#1\right\vert}
\def\set#1{\left\{#1\right\}}
\def\seq#1{\left\langle#1\right\rangle}

\def\sep#1{\left(#1\right)}\def\adf#1{\left[#1\right]}
\def\R{\mathbb R}
\def\N{\mathbb N}
\def\T{\mathbb T}
\def\D{\partial}
\def\eps{\varepsilon}
\def\phi{\varphi}
\def\Re{{\mathrm Re\,}}

\def\defegal{\stackrel{\text{\rm def}}{=}}

\renewcommand{\iint}{ \int \!\!\!\!\! \int}

\newcommand{\dv}{{{\rm d}v}}
\newcommand{\dx}{{{\rm d}x}}

\def\exp{\textrm{e}}
\newcommand{\ddt}{{\frac{{{\rm d}}}{{{\rm d}t}}}}

\newcommand{\dmu}{{\rm d}\mu}
\newcommand{\dnu}{{\rm d}\nu}
\newcommand{\Op}{{\rm Op}}

\newcommand{\Enl}{E_{\textrm{nl}}}
\newcommand{\ds}{{\rm d}s}
\newcommand{\Gl}{G_\ell}
\newcommand{\s}{\sigma}

\def\seta{\seq{\eta}}
\def\sxi{\seq{\xi}}
\def\exp{\textrm{e}}
\def\hf{\widehat{f}}

\author{Fr\'ed\'eric H\'erau}
\address{Laboratoire de Math\'ematiques J. Leray, UMR  6629 du CNRS, Universit\'e de Nantes, 2, rue de la Houssini\`ere, 44322 Nantes Cedex 03, France}
\email{frederic.herau@univ-nantes.fr}

 \title[Introduction to hypocoercive methods]{Introduction to hypocoercive methods and  applications for simple linear inhomogeneous kinetic models. \\ \bigskip Course given at Morningside center of Mathematics in October 2016. }
\thanks{The author wants to thanks professor Ping Zhang and his colleague Wen Deng for their great hospitality at the Morningside Center of Mathematics of Beijing, and Zeinab Karaki for a careful reading of a first version of these notes.}

\begin{document}
\frontmatter

 \begin{abstract} In this lectures given at the Morning side center of Mathematics, Beijing, we present in a very simple framework Hilbertian hypocoercive methods in the case of 1$d$ kinetic inhomogeneous equations, and some illustrations concerning short time or long time behavior in a linear or non-linear perturbative setting.
\end{abstract}

%\thanks{F.H. is supported by the grant "NOSEVOL" ANR-2011-BS01019-01.}
\maketitle
\mainmatter

%\Gr
% \begin{acknowledgements}
%
%\end{acknowledgements}
%\Bk

%
%\Rd preliminary notes (biblio missing, second reading etc...). To be published in a special book of lectures of the morningside center - 2016. confidential. \Bk

\tableofcontents

\section{Introduction}

In this series of lectures, we are interested in inhomogeneous kinetic equations of the following form
  \begin{equation} \label{IHFPF}
 \D_t F + v \D_x F - \D_x V(x) \D_v F = \lll(F), \qquad F|_{t=0} = F_0,
 \end{equation}
 where $F= F(t,x,v)$ is the density of presence of a system of particules, and $t\geq 0$, $d\in \N^*$, $x \in \T^d$ or $\R^d$ and $v \in \R^d$. The unknown function $F(t,.,.)$ is a priori in $L^1 (\dx\dv)$ (space of densities of probability) for each $t\geq 0$. This type of equation modelizes the evolution of a system of particles (plasma, galaxies, ...).

The objects are the following: the (external) potential $V(x)$ is supposed to be such that  $e^{-V} \in L^1(\dx)$ (we will in particular consider the case $V=0$ when $x \in \T^1= \T$). We  pose
$$
\mu(x) = \frac{ e^{-V(x)} }{\int e^{-V(x)}dx}, \qquad \dmu = \mu \dx.
$$
Symmetrically we consider the Gaussian in velocity
$$
\nu(v) = \frac{1}{\sep{2\pi}^{d/2}} \exp^{-v^2/2}, \qquad \dnu = \nu \dx.
$$
The collision kernel $\lll$ is  acting only in velocity.  It has the two fundamental properties:
 \begin{equation} \label{fund1}
 \int \lll(F) dv = 0, \qquad \textrm{ and } \qquad \lll(\nu) = 0.
 \end{equation}

 A direct consequence of the first property is
 $$
 \iint F(t,.) dxdv = \iint F_0 dxdv
 $$
 which corresponds to the conservation of mass of particules.

 As examples of collisions kernels $\lll$, we can mention Fokker-Planck, Boltzmann or Landau operator, but in this pedagogical course we will essentially focus on the fokker-Planck case and the linear Boltzmann case, which are linear and have a one dimensional kernel $\ker{\lll} = span(\nu)$ (in a good functional space, and when seen as operators in velocity only).

 \begin{enumerate}
 \item the Fokker Planck equation writes
 \begin{equation*}
 \D_t F + v \D_x F - \D_x V(x) \D_v F = \D_v(\D_v+v)F , \qquad (original \ \ (FP))
 \end{equation*}
 \item the linear Boltzmann equation reads
  \begin{equation*}
 \D_t F + v \D_x F - \D_x V(x) \D_v F =  \rho \nu -F.\qquad   (original \ \ (BL) )
 \end{equation*}
 where $\rho(x) = \int F(x,v) dv$  and the product $\rho(x) \nu(v) $ is sometimes called the local Maxwellian.
 \end{enumerate}

 This is direct to check that in both cases the two properties \eqref{fund1} are satisfied. Mention just now that these two equations
 seem to have bad properties with respect to standard PDE tools: there is a transport (hyperbolic) part by the vector field
 $$
 X_0 = v \D_x - \D_x V(x) \D_v
 $$
 and a part which is either non-diffusive ((BL) case) or diffusive only in velocity  ((FP) case). We shall see later that $\lll$ has in fact good spectral properties in an adapted subspace of $L^1$, beeing selfadjoint and non-negative. But at this point the two operators seem to be neither elliptic nor selfadjoint.

 One of the central object in the analysis of kinetic equations is the so-called Maxwellian, defined by
 $$
 \mmm(x,v) = \mu(x)\nu(v) = \frac{ e^{-(V(x)+ v^2/2) }}{\int e^{-(V(x)+ v^2/2) }dxdv}.
 $$
 This is direct to see that $\mmm$ is an equilibrium for the two equations (it will appear to be the only one for our models). Indeed by direct computations we have
 $ \D_t \mmm = 0$, $ X_0 \mmm =0 $  since $\mmm$ is a function of the classical hamiltonian  $V(x)+ v^2/2$ of Hamiltonian vector field   $X_0$  and (with a slight abuse of notations) $\lll(\mmm(x,v)) = \mu(x) \lll(\nu(v)) = 0$. A natural related question is then:

\begin{center}\begin{framebox} {
 Question 1 : does $F(t,x,v) \longrightarrow \mmm$ and at which rate ? }
\end{framebox} \end{center}

\bigskip
This question raises many problems : in what space ? which convergence ?  etc... As a first answer we present now a light version of the so-called H-theorem: For a smooth sufficiently decaying  (in space and velocity) and positive density of probability, we denote
$$
H(F,\mmm) \defegal \iint F \cdot \ln \sep{\frac{F}{\mmm}} dxdv = \iint F \ln F  dxdv + \iint \sep{\frac{v^2}{2} + V(x)} F (t,x,v) dxdv .
$$
then we directly check that $H(\mmm,\mmm) = 0$ and that it is the only one, and that
$$
H(F,\mmm) = \iint F \cdot \ln \sep{\frac{F}{\mmm}} dxdv = \iint \sep{  \frac{F}{\mmm}  \ln \sep{\frac{F}{\mmm}} - \sep{\frac{F}{\mmm}} + 1} \mmm dxdv \geq 0
$$
Suppose now that $F(t,\cdot)$ is a density of probability solution of the kinetic equation with the same properties as above. Then we check that
$$
\ddt H(F(t),\mmm) \leq 0.
$$
Indeed e.g. for the Fokker-Planck case
\begin{equation*}
\begin{split}
& \ddt H(F(t),\mmm) \\
 & = \frac{d}{dt} \sep{ \iint F \ln F  dxdv + \iint \sep{\frac{v^2}{2} + V(x)} F (t,x,v) dxdv} \\
& = \iint \D_t F (\ln F + 1) dxdv + \iint \D_t F \sep{\frac{v^2}{2} + V(x)}  dxdv \\
& - \iint  X_0 F (\ln F + 1) dxdv - \iint X_0 F \sep{\frac{v^2}{2} + V(x)}  dxdv \\
 &  \qquad + \iint \D_v(\D_v +v) F (\ln F + 1) dxdv + \iint \D_v(\D_v+v)F \sep{\frac{v^2}{2} + V(x)}  dxdv
\end{split}
\end{equation*}
The first two terms are $0$ since  $X_0 \sep{\frac{v^2}{2} + V(x)} = 0$ and $X_0 F (\ln F + 1)= X_0( F\ln F)$. For the last two ones, we get by IPP again
\begin{equation*}
\begin{split}
\ddt H(F(t),\mmm)
 &  = - \iint (\D_v +v)F \D_v F \frac{1}{F}  dxdv + \iint (\D_v+v)F v  dxdv \\
 & = - \iint \frac{|(\D_v +v)F|^2}{F}  dxdv \defegal - D(F(t), \mmm) \leq 0
\end{split}
\end{equation*}
We could do the same for the linear Boltzmann equation and more generally for all standard kinetic models. The term $D(F,\mmm)$ is called the dissipation of entropy term.

This result answers partially  to the question, at least in the space corresponding to finite entropy densities. Now this gives no idea of a possible rate, since there is no hope to have $D(F,\mmm) \geq H(F,\mmm)$ as in the homogeneous case :  the gain in variable $x$ is missing. This is the place where hypocoercivity methods enter, and we shall later completely prove the exponential time decay in an adapted space we now build thanks to perturbation considerations.

\bigskip The kinetic equations in which we are interested in this course are linear, but in a perturbative context it is useful to consider a linearization near the maxwellian $F$, and we pose
$$
F = \mmm + \mmm f
$$
where $f$ is now supposed to be small.  A first fundamental remark is that owing to the conservation of mass we have
$$
\iint f d\mu d\nu = \iint f\mmm dxdv = \iint F dxdv - \iint \mmm dxdv = 1-1 = 0
$$
Note this can be understood as saying that $f \perp 1$ in $L^2 (d\mu d\nu)$ and this is perhaps the first time this space appears in the perturbative analysis.
For the (BL)  and (FP) equations this approach is essentially transparent since the  equations are linear. We just have to deal with the multiplication by $\mmm$ and we get equations of the form $\D_t f + v\D_x f = L(f)$.
%\begin{equation} \label{eq:rescaled}
%\D_t f + v\D_x f = L(f)$.
%\end{equation}
Precisely,

\begin{enumerate}
 \item the Fokker Planck equation writes
 \begin{equation*}
 \D_t f + v \D_x f - \D_x V(x) \D_v f = -(-\D_v+v) \D_v f , \qquad  (FP)
 \end{equation*}
 \item the linear Boltzmann equation reads
  \begin{equation*}
 \D_t f + v \D_x f - \D_x V(x) \D_v f = r  -f .\qquad  (BL)
 \end{equation*}
 where $r(x) = \sep{\int f(x,v) d\nu}$.
 \end{enumerate}

Note that we strongly used that $X_0 \mmm = 0$ in the computation. These equations have a priori essentially the same form. Now for the relative entropy the change is very enligthening. We do below only formal computations, and use that $f$ is "small" :
\begin{equation*}
\begin{split}
H(F,\mmm) & = \iint \mmm( 1+f) \ln( 1+f) dxdv \\
& \sim \iint \mmm( 1+f) (f - \frac{f^2}{2} + o(f^2)) dxdv \\
& \sim \iint f d\mu d\nu + \frac{1}{2} \iint f^2 d\mu d\nu + \iint o(f^2) d\mu d\nu \\
& \sim \frac{1}{2} \iint f^2 d\mu d\nu
\end{split}
\end{equation*}
In this perturbative context, the  relative entropy is therefore the $L^2(d\mu d\nu)$ norm  and we shall in this course study indeed the decay of the $L^2(d\mu d\nu)$ norm of the solution $f$ of the (perturbative) (FP) and (BL)  equations. This is a small space with respect to $L^1$ but note anyway that the embedding
 $$
 f \in L^2(\dmu\dnu) \hookrightarrow L^1 (dxdv) \ni F
 $$
  is of norm 1. Note also that either using this formal perturbative strategy or by doing exact computations, the entropy-entropy-dissipation inequality reads
\begin{equation*}
\begin{split}
\ddt \frac{1}{2} \norm{f(t)}^2
 & = - \norm{ \D_v f}^2 \leq 0.
\end{split}
\end{equation*}
This is again not sufficient  to get some explicit decay rate  and we shall build later modified entropies, in this Hilbert context and prove completely (in simple cases)  that there is indeed an exponential time decay rate   towards the equilibrium.

 Before stating stating the result in our simple context, we mention a series of (non optimal) hypotheses  on the potential $V$. The first one is that $V$ is smooth and with derivatives of order 2 or more bounded. The second additional hypothesis is that there is a Poincare inequality for $d\mu$ : there exists a constant $c_P$ such that for all
 $$
 \forall \phi \in H^1(d\mu), \qquad c_P \norm{\phi - \seq{\phi}}^2 \leq \norm{\D_x \phi}^2
 $$
 where $\seq{\phi} = \int \phi d\mu$ and $H^1(d\mu)$ is the space of functions whose differential are in $L^2(d\mu)$. Anyway in this lectures we shall essentially focus on the case when $V=0$, $d=1$ and $(x,v) \in \T\times \R$, and note that the preceding hypotheses are satisfied in this case. The following theorem concerns the return to the equilibrium and will be proven in the next section.

 \begin{thm} \label{hypocoercivite} Suppose that $V = 0$, $d=1$, $(x,v) \in \T\times \R$. Let $f$ is a solution of  (BL)  or  (FP) with $ \seq{ f_0 } = 0$. Then there exists $\kappa >0$ and $C>0$ explicit (independant of $f_0$) such that for all $t\geq 0$, $\norm{f(t)} \leq c e^{-\kappa t}\norm{f_0}$.
 \end{thm}

(We point out that in the preceding statement, we use the notation
$ \seq{f_0} = \iint f_0 \dmu\dnu$ for the mean in space and velocity variables.)

We will present two alternative proofs of the exponential time decay, first in $H^1(\dmu\dnu)$ for  (FP) (see the precise statement in Section 2.2) and then in $L^2$ (Section 2.4) for  (FP) and  (BL). This type of result was first proven in $L^2$ (with the explicit constants and a general external potential $V$) in \cite{HN04} but the method here was first developed in an $L^2$ context in \cite{Her06} (see also \cite{DMS15}), then in an  $H^1$ context in \cite{MN06} and  in a general framework in \cite{Vil09}. The presentation here of the $H^1$ version come from \cite{DHL16} where the discrete case is also analyzed.
The presentation of the $L^2$ proof here follows a so-called \it micro-macro \rm scheme; it is essentially new and follows in a simpler situation the proof proposed in \cite{DHMS16} for more complex model of Boltzmann type.

\bigskip
Of course we first need to precise what we mean by solution of the equation:  here this means that $f(t) = e^{-tP} f_0$ where $P = X_0-L$. This supposes that we have been able to apply Hille Yosida Theorem to operator $P$ with an appropriate domain $D(P)$, and we shall assume this in the following. In particular shall always have
$$
f \in \ccc^1(\R, L^2(\dmu \dnu)) \cap \ccc^0(\R, D(P)).
$$
Note that for  (FP) this is not an easy result (see e.g. \cite{HelN04} in a general context). We shall also assume the same in $H^1(\dmu\dnu)$.

We mention here the fundamental properties of the (rescaled) collision kernel $L$. The first one is that it is selfadjoint in $L^2(\dnu)$ and the second main additional property is that it has a spectral gap: there exists a constant $c_L$ such that
$$
 \forall \phi \in D(L), \qquad c_L \norm{\phi - \seq{\phi}}^2 \leq - \int L\phi(v) \phi(v) \dnu,
 $$
where $D(L)$ is the domain of operator $L$ in $L^2(d\nu)$ and where we denote again $\seq{\phi} = \int \phi d\nu$. This property is trivially satisfied in the (BL)  case since in this case $-\int L\phi  \phi \dnu= \norm{\phi-\seq{\phi}}^2$. In the (FP) case we have
$$
-\int L\phi \phi \dnu= \norm{\D_v \phi}^2 \geq \norm{\phi}^2,
$$
where the last inequality is the usual Poincar\'e inequality in velocity for the measure $d\nu$. The spectral gap is indeed equal to $1$.

\bigskip
Now a second question is natural in the context of kinetic equations when the collision kernel has in addition diffusive properties.
 We focus now only on the (FP) case. As we mentioned before, the involved operator has elliptic properties only in velocity and not in the space variable. One can nevertheless ask about the second main question of  this course :
 \begin{center}
 \begin{framebox}
 { Question 2 : Is $F(t)$ when $t>0$ more regular that $F_0$ and how much ?}
 \end{framebox}
 \end{center}

 \bigskip
It appears that there are indeed regularizing properties for  (FP), and this  is one manifestation of the so-called hypoellipticity of operator $P=X_0-L$. We now give a result in the simple situation $V=0$, $d=1$, $(x,v) \in \T\times \R$, and we with $f$ where  $F=\mmm + \mmm f$.

\begin{thm} \label{hypoellipticite} Suppose that $V=0$, $d=1$ and $(x,v) \in \T\times \R$ with $f_0 \in L^2(\dmu\dnu)$. Then the solution $f$ of  (FP) satisfies $f(t) \in H^1$ for all $t>0$ and precisely there exists $C>0$ explicit (independant of $f_0$) such that for all $t\in (0,1]$,
$$
\norm{\D_v f(t)} \leq \frac{C}{ t^{1/2}} \norm{f_0}, \qquad \textrm{and} \qquad \norm{\D_x f(t)} \leq \frac{C}{ t^{3/2}}\norm{f_0}.
$$
\end{thm}

The result and the method were first presented in \cite{Her07} in a more general context. In fact in can be proven that $f(t) $ is in the Schwartz space for more general intial data using alternative Cauchy contour methods (see \cite{HN04}) but what is interesting here is the rate and the simplicity of the method that can be adapted to various models and situations. The question of the regularization properties of kinetic equations has a long history coming back to hypoellipticity results by Kohn \cite{Koh73}
or H\"ormander \cite{Hor67}. In the result above we again focus on a very simple case $V=0$, $d=1$ and $(x,v) \in \T\times \R$ that could be treated with explicit formulas (see e.g. \cite{Hor95}). The result here were inspired by studies of subelliptic semi-groups in the spirit of  \cite{CSV92} or \cite{FS86}. the strategy presented here has been adapted to a very large family of kinetic diffusive equations in various functional contexts (see Theorem \ref{FKEthm} below and e.g.  \cite{GMM13} in large space, \cite{Vil09} in $L \log L$ spaces or \cite{PZ16} in discrete cases). The homogeneous fractional Fokker-Planck case has been studied in \cite{Tri15} and the general Boltzmann without cutoff case in \cite{HTT16b}.

\bigskip
As an illustration  we give now a result for the fractional Kolmogorov equation  which has a partially diffusive collision kernel. Let $s\in (0,1]$, the fractional Kolomogorov equation  writes
\begin{equation*}
 \D_t f + v \D_x f  = -( 1 -\Delta_v)^{s}f , \qquad (FK)
 \end{equation*}
For convenience we introduce
the following strictly positive operators
$$
\Lambda_v^2 =   1 -\Delta_v,\qquad \Lambda_x^2 =  1 -\Delta_x
$$
 and the associated family of Fourier multipliers
 $$
\Lambda_x^\alpha = ( 1 -\Delta_x)^{\alpha/2}, \qquad  \Lambda_v^\beta =  ( 1 -\Delta_v)^{\beta/2}, \qquad \qquad  \alpha, \beta \in \R
$$
 which act on a function in $\sss(\T\times \R)$ or $\sss(\R\times \R)$ in the following way
 $$
 \widehat{ \Lambda_x^\alpha f} (\xi, \eta) = (1+\xi^2)^{\alpha/2} \widehat{f}(\xi, \eta), \qquad \widehat{ \Lambda_v^\beta f} (\xi, \eta) = (1+\eta^2)^{\beta/2} \widehat{f}(\xi, \eta)
 $$
 where the hat corresponds to the Fourier transform in both $x$ and $v$ variables. For convenience we will denote $\sxi = (1+\xi^2)^{1/2}$ and  $\seta = (1+\eta^2)^{1/2}$. We also introduce the corresponding Sobolev spaces
 $$
 H^{\alpha,\beta} = \set{ f \in \sss', \quad \Lambda_x^\alpha \Lambda_v^\beta f \in L^2},
 $$
 and we denote by $\norm{\cdot}_{\alpha,\beta}$ the corresponding norm. For $r\in \R$ We also denote
 $$
 H^r = H^{0,r} \cap H^{r,0}
 $$
 the isotropic Sobolev space.
 With the notations introduced above, the fractional Kolmogorov equation reads
 \begin{equation} \label{FKE}
 \D_{t}f+v\D_{x}f + \Lambda_v^{2s}f =0
 \end{equation}
 and a natural question is wether $f$ benefits from some regularization
 induced by the  elliptic  properties of $\Lambda_v^{2s}$.
 The main result concerning the fractional Kolmogorov equation is the
 following:

 \begin{thm} \label{FKEthm}
 Let $r \in \R$ and $f$ be a solution of (FK)  with initial data $f_0 \in H^{r,0}_{x,v}$. Then, there exists a constant $C_r >0$ independent of $f_0$ such that for all $t \in (0,1]$,  we have
 $$
 \norm{f(t)}_{r,s} \leq \frac{C_r}{t^{1/2}} \norm{f_0}_{r,0}
\textrm{  and  }
 \norm{f(t)}_{r+s,0} \leq \frac{C_r}{t^{1/2+s}} \norm{f_0}_{r,0}.
 $$
 \end{thm}

Note first that when $r=0$ and $s=1$ we essentially recover Theorem \ref{hypoellipticite}. This version for the fractional Kolmogorov $s \in (0,1]$ comes from \cite{HTT16a} here stated in a $1d$. This type of result is of great use in the proof of the return to equilibrium in large functional spaces of solutions of inhomogeneous kinetic equations following the general method  presented in  \cite{GMM13}. The homogeneous fractional Fokker-Planck case has been studied in \cite{Tri15} where regularization properties in velocity are investigated thanks to a fractional Nash inequality. The present result is one of the stone of the proof in the inhomogeneous Boltzmann without cutoff case proposed in \cite{HTT16b}.
 In the present lecture we pay attention to give a proof \it not \rm using
 any kind of pseudodifferential tool (only Fourier multiplier), although the proof is deeply of microlocal inspiration.

\bigskip

The proof of the three preceding results (Theorems \ref{hypocoercivite}, \ref{hypoellipticite} and \ref{FKEthm}) have remarkable similarities in spirit and in
the computations. The name \it hypocoercivity \rm (hidden coercivity of operator $P$) is the brother of the name \it hypoellipticity \rm (hidden regularization effect or hidden ellipticity) although the results have a completely different nature: the first  one is of spectral inspiration and the other one concern diffusion properties. This is remarkable that similar methods
yield results of such a different nature.

\bigskip
As an application to the previous hypocoercive methods, we propose in the last part of these lectures the study of the Cauchy problem and the trend to the equilibrium for a mollified Vlasov-Poisson-Fokker-Planck model in a perturbative situation. This non-linear model of mean-field type reads
\begin{equation*}
\left\{
\begin{array}{l}
 \D_t F + v \D_x F +\Enl(F) \D_v (F-\mmm) = \D_v(\D_v+v)F,\\
% \pm \Enl(f) = -\D_x \Vnl (x), \\
%  \textrm{ with } -\Delta_x \Vnl(x) = \rho - 1 \textrm{ and } \rho(t,x) = \int F(t,x,v) \dv \\
  \pm \Enl(F) = K*(\rho - 1) \qquad \textrm{ and } \rho(t,x) = \int F(t,x,v) \dv \\
 F|_{t=0} = F_0, \qquad \iint F_0 \dx\dv = 1
 \end{array}
 \right. \qquad (MVPFP)
 \end{equation*}
 where $K$ is supposed to be smooth and in $ L^\infty(\dx)$.  Using both short and long time hypocoercive estimates, we get the following result stated here again in the very simple situation  $d=1$ and $(x,v) \in \T\times \R$, and therefore $\mmm(x,v) = \nu(v)$. The definition of mild solution is given is Section 4 below.

\begin{thm} \label{MVPFPbis} Suppose $d=1$,  and work on $\T$ in the space variable. Then there exists $\eps_0$, $C_0$ and $\kappa_0$ explicitely computable such that the following happens; denote $F = \mmm + \mmm f$ and suppose $\seq{f_0}=0$:
 \begin{enumerate}
 \item If $\norm{f_0} \leq \eps_0$ then there exists a unique mild solution  $(f,E) \in \ccc^0(\R_+, L^2(\dmu\dnu) ) \times \ccc^0(\R^+, L^\infty(\dx))$ to the mollified Vlasov-Poisson-Fokker-Planck system (MVPFP).
     \item We have in addition for all $t\geq 0$,
     $$
     \norm{f(t)}_{L^2(\dmu\dnu)} \leq  C_0 \exp^{-\kappa_0 t} \quad \textrm{ and } \quad \norm{E(f(t))}_{L^\infty(\dx)} \leq  C_0 \exp^{-\kappa_0 t}
     $$
     \end{enumerate}
     \end{thm}

This result is a simplified version of the one presented in \cite{HT14} where a  result with more general external potentials and singular mean-field potentials of Coulomb type  is proposed in any dimension. This type of problem has a very long history and we only mention the two major non-perturbative results (in the case of no external potential) concerning the Coulomb case:
the $2d$ case was treated in \cite{Deg86} and and the $3d$ case was treated in \cite{Bou93} with the help of the explicit heat kernel of the linear Fokker-Planck equation without potential. Mention that as soon as an external potential is present, no explicit formula for the heat kernel is available. The interest of the approach presented here is that it is not based on any explicit formula.

\bigskip

The plan of the course is the following. For simplicity and pedagogy, we  focus on the simple case $d=1$, $V=0$ and $(x,v) \in \T \times \R$ or $\R\times \R$. In Section 2 we will first prove an enlightening version of the exponential trend to the equilibrium in $H^1$ (Corollary \ref{hun}) for the (FP) case,  and then give the proof in $L^2 (\dmu \dnu)$  for both models simultaneously (Theorem \ref{hypocoercivite}).  In section 3 we will study the short time regularization and propose the proof of Theorem \ref{hypoellipticite} and Theorem \ref{FKEthm}. In a last section, we give the application to the simple nonlinear (MVPFP) problem and prove Theorem \ref{MVPFPbis}

\section{Trend to the equilibrium}

In this section we  introduce the notion of hypocoercivity, with two proofs of the trend to the equilibrium (Corollary \ref{hun} in $H^1(\dmu\dnu)$ and  Theorem \ref{hypocoercivite} in $L^2(\dmu\dnu)$), both enlightening the strategy.  In all what follows we focus on the perturbed equations near the equilibrium with unknown function $f$ with mean $0$ and with Hilbertian entropies, as presented in the introduction.

\subsection{The coercive case}
Although nearly trivial it may be interesting to focus on the coercive case first.
The equation of evolution is only in velocity $v \in \R^d$ and time $t\geq 0$ and it reads in this case
\begin{equation} \label{coer}
\left\{
\begin{array}{l}
\D_t f - L f = 0 \\
f|_{t=0} = f_0, \qquad \seq{f_0} = 0,
\end{array}
\right.
\end{equation}
where $\seq{f} = \int f \dnu$ and
$$
Lf= -(-\D_v + v)\D_v f\qquad  (FP) \qquad \textrm{or} \qquad L f=-( f-\seq{f}) = -f \qquad  (BL)
$$
The entropy is then
$$
\hhh(f) = \frac{1}{2} \norm{f}^2
$$
and if $f$ is the solution of \eqref{coer} (given by Hille Yosida Theorem), we get
\begin{equation} \label{zeinab}
\ddt \hhh(f) = \seq{Lf,f} \leq \left\{\begin{array}{l} -\norm{f}^2 \ \textrm{ or } \\ -\norm{\D_v f}^2 \end{array} \right. \leq -\norm{f}^2 = - \hhh(f)
\end{equation}
where this is direct for  (BL)  and we used Poincar\'e inequality for  (FP).
This immediately gives by
$$
\norm{f(t)}^2 \leq \exp^{-t} \norm{f_0}^2
$$ and the result is proven.

\subsection{Hypocoercivity for  (FP) in $H^1$}
In the inhomogeneous case, we focus for convenience on the case when $d=1$, $(x,v) \in \T \times \R$ and $V=0$.  We look first at the   (FP) model.
The equation reads
\begin{equation} \label{eq:rescaledFP}
\left\{
\begin{array}{l}
\D_t f + v\D_x f +(-\D_v+v) \D_v f = 0  \\
f|_{t=0} = f_0, \qquad \seq{f_0} = 0
\end{array}
\right.
\end{equation}
where here $\displaystyle \seq{f} = \iint f \dmu\dnu $ (be careful the integration is in both variables).
For convenience we introduce the macroscopic density
$$
r(t,x)= \int f(t,x,v) \dnu \qquad (\textrm{ note that } \seq{r} = 0).
$$
Of course the coercive method does not apply since
the same method only yields
$$
\ddt \frac{1}{2} \norm{ f}^2  \leq -\norm{ f- r}^2,
$$
and we do not recover the full entropy with a minus sign on the right. In order to close such an inequality we introduce the following modified entropy
 for $C>D>E>1$ to be fixed later :
\begin{equation}
  \label{eq:entropyfunc}
  \eee(f) = C\norm{f}^2+D\norm{\partial_vf}^2+E\seq{\partial_v f,\partial_x f}+\norm{\partial_xf}^2
\end{equation}
We will show that for well chosen $C,D,E$, $t\mapsto \eee(f(t))$ is
nonincreasing when $f$ solves equation (FP) with initial datum $f_0\in
H^1(\dmu\dnu)$.
We first prove that $\eee$ is equivalent to the $H^1(\dmu\dnu)$-norm.
\begin{lem}\label{lem:equiv1}
  If $E^2<D$ then
$$%  \begin{equation}
%    \label{eq:equiv}
    \dfrac{1}{2}\norm{f}_{H^1}^2\leq\eee(f)\leq 2C\norm{f}_{H^1}^2
$$%  \end{equation}
\end{lem}

\preuve
  We use the Cauchy-Schwarz inequality and observe that
  \begin{equation*}
    2\abs{E\seq{\partial_v f,\partial_x f}}\leq E^2\norm{\partial_v f}^2+\norm{\partial_x f}^2
  \end{equation*}
which implies for all $f\in H^1(\dmu\dnu)$
\begin{multline*}
  %\underbrace{C}_{1/2\leq}
  C \norm{f}^2
  +%\underbrace{(D-E^2/2)}_{1/2\leq E^2/2\leq}
  (D-E^2/2)\norm{\partial_vf}^2
  +\sep{1-1/2}\norm{\partial_xf}^2\leq \hhh(f) \\ \leq
  C\norm{f}^2+
  %\underbrace{(D+E^2/2)}_{\leq D+D/2\leq 3C/2\leq    2C}
    (D+E^2/2)\norm{\partial_vf}^2+%\underbrace{3/2}_{\leq 3C/2\leq 2C}
    \frac{3}{2}\norm{\partial_xf}^2.
\end{multline*}
This implies the result if $E^2<D$.
\fin

\begin{prop}\label{prop:decrexp1}
   There exists $C,D,E$ and $\kappa>0$ such that for all $f_0\in
  H^1(\dmu\dnu)$ with $\seq{f_0}=0$ the solution of \eqref{eq:rescaledFP} satisfies
  \begin{equation*}
   % \label{eq:estimH}
    \forall t>0,\qquad \eee(f(t))\leq \eee(f_0)e^{-\kappa t}
  \end{equation*}
\end{prop}

\preuve
  We compute separately the time derivatives of the four terms defining
  $\eee(f(t))$. Omitting the dependence on $t$, the first one reads
  \begin{align*}
    \dfrac{d}{dt}\norm{f}^2&=2\seq{\partial_tf,f}
=-2\underbrace{\seq{v\partial_xf,f}}_{=0}-2\seq{(-\partial_v+v)\partial_vf,f}
=-2\norm{\partial_vf}^2
  \end{align*}
The second term writes
\begin{align*}
  \dfrac{d}{dt}\norm{\partial_vf}^2 &=2\seq{\partial_v(\partial_tf),\partial_vf}\\
&=-2\seq{\partial_v(v\partial_xf+(-\partial_v+v)\partial_vf),\partial_v f}\\
&=-2\underbrace{\seq{v\partial_x\partial_vf,\partial_vf}}_{=0}-2\seq{\adf{\partial_v,v\partial_x}f,\partial_vf}-2\seq{\partial_v(-\partial_v+v)\partial_vf,\partial_vf}.
\end{align*}
We again use the fact that $v\partial_x$ is a skewadjoint operator  and the fundamental relation
$\adf{\partial_v,v\partial_x}=\partial_x$ and we get
\begin{equation*}
  \dfrac{d}{dt}\norm{\partial_vf}^2 =-2 \seq{\partial_xf,\partial_vf}-2\norm{(-\partial_v+v)\partial_vf}^2.
\end{equation*}
The time derivative of the third term can be computed as follows
\begin{align*}
   \dfrac{d}{dt}\seq{\partial_xf,\partial_vf}&=-\seq{\partial_x(v\partial_xf+(-\partial_v+v)\partial_vf),\partial_vf}-\seq{\partial_xf,\partial_v(v\partial_xf+(-\partial_v+v)\partial_vf)}\\
&=-\seq{v\partial_x(\partial_xf),\partial_vf}-\seq{\partial_x\partial_vf,
  \partial_v^2f}-\seq{\partial_xf,\adf{\partial_v,v\partial_x}f}-\seq{\partial_xf,v\partial_x\partial_vf}\\
&\qquad -\seq{\partial_xf,\adf{\partial_v,(-\partial_v+v)}\partial_vf}+\seq{(-\partial_v+v)\partial_vf,\partial_x\partial_vf}.
\end{align*}
Using the fact that $v\partial_x$ is skewadjoint, we have
\begin{equation*}
  \seq{v\partial_x\partial_xf,\partial_vf}+\seq{\partial_xf,v\partial_x\partial_vf}=0.
\end{equation*}
Since $\adf{\partial_v,v\partial_x}=\partial_x$ and
$\adf{\partial_v,(-\partial_v+v)}=1$, we get
\begin{equation*}
   \dfrac{d}{dt}\seq{\partial_xf,\partial_vf}=-\norm{\partial_xf}^2+2\seq{(-\partial_v+v)\partial_vf,\partial_x\partial_vf}-\seq{\partial_xf,\partial_vf}.
\end{equation*}
Finally, observing that $\partial_x f$ also solves \eqref{eq:rescaledFP}, we
obtain for the last term on $\eee(f(t))$ the same estimate as the one we
obtained for the first term:
\begin{equation*}
    \dfrac{d}{dt}\norm{\partial_xf}^2 =-2\norm{\partial_v\partial_xf} ^2.
\end{equation*}
Eventually, we obtain
\begin{align*}
   \dfrac{d}{dt}
  \eee(f)&=-2C\norm{\partial_vf}^2-2D\norm{(-\partial_v+v)\partial_vf}^2-E\norm{\partial_xf}^2-2\norm{\partial_x\partial_vf}^2\\
&\qquad -(2D+E)\seq{\partial_xf,\partial_vf}+2E\seq{(-\partial_v+v)\partial_vf,\partial_x\partial_vf}.
\end{align*}
Only two of the terms do not have a sign. Using the Cauchy-Schwarz inequality,
we observe that
\begin{equation*}
\left\{ \begin{array}{l}  \abs{(2D+E)\seq{\partial_xf,\partial_vf}}\leq \dfrac{1}{2}\norm{\partial_xf}^2+\dfrac{1}{2}(2D+E)^2\norm{\partial_vf}^2 \\
  \abs{2E\seq{(-\partial_v+v)\partial_vf,\partial_x\partial_vf}}\leq
  \norm{\partial_x\partial_vf}^2 +E^2 \norm{(-\partial_v+v)\partial_vf}^2.
  \end{array}
  \right.
\end{equation*}
Therefore, assuming again that $1<E<D<C$, $E^2<D$ and $\dfrac{1}{2}(2D+E)^2<C$, we get
\begin{equation*}
   \dfrac{d}{dt}\eee(f)\leq
   -C\norm{\partial_vf}^2-(E-1/2)\norm{\partial_xf}^2\leq -\dfrac{E}{2}(\norm{\partial_vf}^2+\norm{\partial_xf}^2).
\end{equation*}
Using the Poincar\'e inequality in space and velocity, we
derive
\begin{align*}
  -\dfrac{E}{2}(\norm{\partial_vf}^2+\norm{\partial_xf}^2)&\leq
   -\dfrac{E}{4}(\norm{\partial_vf}^2+\norm{\partial_xf}^2)-\dfrac{E}{4} c_p\norm{f}^2\leq-\dfrac{E}{4}\dfrac{c_p}{2C}\eee(f).
\end{align*}
\fin

%\bigskip
%\preuve[of Theorem ***]  The proof  is direct : Choose $C>D>E>1$ as in Theorem \ref{th:decrexp}. Set $\kappa=E/(4C)$ and $c=\sqrt{2C}$ and
%  apply Theorem \ref{th:decrexp} and Proposition \ref{prop:equiv}. \fin

\begin{cor} \label{hun}
  There exists $c,\kappa>0$ such that for all $f_0\in
  H^1(\dmu\dnu)$ with $\seq{f_0}=0$, the solution of \eqref{eq:rescaledFP}
  satisfies
  \begin{equation*}
    % \label{eq:expest}
    \forall t\geq 0,\quad\norm{f(t)}_{H^1}\leq ce^{-\kappa t}.
  \end{equation*}
\end{cor}
\preuve
  Choose $C>D>E>1$ as in Proposition \ref{prop:decrexp1}. Set $k=E/(4C)$ and $c=\sqrt{2C}$ and
  apply Proposition \ref{prop:decrexp1} and the equivalence of norms in Lemma \ref{lem:equiv1}.
\fin

\remark As a general remark about this proof, mention that main ingredient is the fact that
$\adf{\partial_v,v\partial_x}=\partial_x$ which allows to recover the missing $\D_x$ derivative in the computations. This is of microlocal inspiration: the family of vector fileds $ \D_v$,  $\D_x$ and their first commutators $\D_x$ span the full tangent space at each point $(x,v)$. Note anyway that it was used here to recover a \it spectral \rm result. This is the core the hypocoercivity method.

\subsection{Preliminaires  in the hypocoercive $L^2$ case}

The preceding method seems to apply only in the Fokker-Planck case. In fact it can be adapted in many situations for kinetic models not having any diffusion property. We propose in the next Section a short proof of Theorem \ref{hypocoercivite} for both models (FP) and (BL), coming from the same fundamental idea and in the $\L^2$ case. Again we assume that we are in the simple situation $d=1$, $V=0$ and $(x,v) \in \T\times \R$.

For this we first write the so-called macroscopic equations.
For the following we denote
\begin{eqnarray*} %\label{macro}
&r(x) = \int f(x,v) \dnu \ \textrm{ (local mass), } \\
& m(x) =  \int v f(x,v) \dnu \ \textrm{ (local moment). }
\end{eqnarray*}
Note that $r$ and $mv$ can be interpreted as orthogonal projections of $f$. The equation reads
\begin{equation} \label{eq:rescaled}
\left\{
\begin{array}{l}
\D_t f + v\D_x f = L(f) \\
f_{t=0} = f_0, \qquad \seq{f_0} = 0
\end{array}
\right.
\end{equation}
with $L$ given by
$$
L(f) = -(-\D_v+ v) \D_v \quad  ((FP) \ case) \qquad \textrm{ or } \qquad Lf = r  - f \quad ((BL) \ case).
$$
The functions $r$, $m$  are elements of $L^2(\dmu)$ and we note that
the functions $r(x)$ and $m(x) v$ are orthogonal in $L^2(\dmu\dnu)$ just because $1$ and $v$ are orthogonal (and normalized) with respect to the scalar product on $L^2(\dnu)$. Using this fact again we make a so-called micro-macro decomposition of
any function $f \in L^2(\dmu\dnu)$ and we denote
$$
f(x,v)= r(x) +  h(x,v).
$$
Now considering a solution $f$ of the equation \eqref{eq:rescaled} we get the so-called macroscopic equations

\begin{lem}
let $f=r+h$ be the solution. Then we have
$$
\D_t r = \Op_1(h), \qquad \D_t m = -\D_x r + \Op_1(h),
$$
where ${\Op_1(g)}$ denotes a generic bounded operator form $L^2$ to $H^{1,0}$.
\end{lem}

\preuve
For the first one we just sum in velocity so that
$$
\D_t \int f \dnu + \D_x \int vf \dnu = \int Lf \dnu = \seq{Lf, 1} = \seq{f, L1} = 0
$$
where we used that $L1 = 0$ and that $L$ is selfadjoint. For the second one we first multiply by $v$ before doing the integration and we get
$$
\D_t \int vf \dnu + \D_x \int v^2 f \dnu = \int Lf v\dnu.
$$
Now we use that
$$
\int v^2 f \dnu = \int (v^2-1) f \dnu + \int f\dnu = r + \int (v^2-1) h \dnu
$$
and that (for both  (FP) and  (BL)  collision kernels)
$$
\int Lf v\dnu = \seq{Lf, v} = \seq{f, Lv} = \seq{f, v} = \int f v\dnu = \int h v \dnu
$$
and we get the result. \fin

The fact that remainder term are of type $\Op_1(h)$ is not comfortable for the following and we now introduce an operator allowing to go back in $L^2$.
Recall the notation
$$
\Lambda_x^2 = -\Delta_x +1.
$$
Then by standard elliptic results, $\Lambda_x^2$ is elliptic in $L^2(\dmu)$, selfadjoint, invertible from $H^2(\dmu)$ to $L^2(\dmu)$, with $\Lambda_x^2 \geq Id$ and we have the following spectral results
\begin{lem}
For all $\phi \in L^2(\dmu)$ we have
$$
\norm{\Lambda_x^{-1} \D_x \phi} \leq \norm{\phi}, \qquad \norm{\Lambda_x^{-2} \D_x \phi} \leq \norm{\phi}
$$
 and we have the following $L^2(\dmu)$-Poincar\'e inequality :
 $$
\frac{c_p}{1+c_p} \norm{\phi-\seq{\phi}}^2 \leq  \norm{\Lambda_x^{-1} \D_x \phi}^2
$$
\end{lem}

\preuve For the first two inequalities we just observe that for $\psi$ smooth we have
$$
\norm{\D_x \psi}^2  \leq \seq{ (-\Delta_x + 1) \psi , \psi} \leq \norm{ \Lambda \psi}^2.
$$
 Then taking $\phi = \Lambda_x \psi$ gives that operator $\D_x \Lambda_x^{-1}$ is bounded by $1$. The result comes after taking the adjoint operator and extending the result into $L^2(\dmu)$.

For the Poincare inequality, take $\phi \in L^2(\dmu)$ with $\seq{\phi} = 0$ and  notice that
$$
\norm{\Lambda_x^{-1} \D_x \phi}^2 = \seq{ (-\Delta_x + 1)^{-1} (-\Delta_x) \phi, \phi}.
$$
We use now  the spectral theorem with the eigenfunctions of $-\Delta_x$ except $1$ (which is related to the eigenvalue $0$) : denoting $0 = \lambda_1 < \lambda_2 < \lambda_3 <  ...$ (where  $\lambda_2 = c_p$) the increasing sequence of eigenvalues of $-\Delta_x$  we get
$$
\norm{\Lambda_x^{-1} \D_x \phi}^2 \geq \min_{k \geq 1} \frac{\lambda_k}{ \lambda_k +1} \norm{\phi}^2 \geq \frac{\lambda_2}{ \lambda_2 +1} \norm{\phi}^2
$$
since the function $s \longmapsto s/(1+s)$ is increasing over $\R^+$. This gives the result. \fin

\remark Note that in the preceding lemma we could have used the $L^2(\dmu\dnu)$ scalar product instead of the $L^2(\dmu)$ one. Indeed the embedding $L^2(\dmu) \longrightarrow L^2(\dmu\dnu)$ is of norm $1$ since $\dnu$ is a probability measure, and for functions depending only of the $x$ variable, we have clearly
$$
\norm{\phi}_{L^2(\dmu)}^2 = \int \phi^2(x) \dmu = \iint \phi^2(x) \dmu\dnu = \norm{\phi}_{L^2(\dmu\dnu)}^2.
$$
The same is true for scalar products and also concerning the $v$ variable. In the following we shall write them all $\seq{.,.}$ since they coincide when having a meaning.
We shall also denote by $\norm{\cdot}$ the $L^2$ norm.

\subsection{Hypocoercivity in $L^2$}

We are now in position to build the new entropy.
For $\eps >0$ to be fixed later, we denote
$$
\fff(f) = \norm{f}^2 + \eps \seq{\Lambda_x^{-2} \D_x r, m}
$$

We first prove that $\fff$ is equivalent to the $L^2(\dmu\dnu)$-norm.
\begin{lem}\label{lem:equiv2}
  If $\eps \leq 1/2$ then
  \begin{equation*}
   % \label{eq:equiv}
    \dfrac{1}{2}\norm{f}^2\leq\fff(f)\leq 2\norm{f}^2
  \end{equation*}
\end{lem}

\preuve
  We use the Cauchy-Schwarz inequality and use the preceding lemma. We get
  \begin{equation*}
   \abs{\seq{\Lambda_x^{-2} \D_x r, m}} \leq \norm{\Lambda_x^{-2} \D_x r}\norm{m} \leq \norm{r}\norm{m} \leq \norm{f}^2.
  \end{equation*}
This gives the result if $\eps \leq 1/2$.
\fin

Now we can prove the main hypocoercive result leading to Theorem \ref{hypocoercivite}.

\begin{prop}\label{prop:decrexp2}
   There exists $\kappa>0$ such that for all $f_0\in L^2(\dmu\dnu)$ with  $\seq{f_0}=0$, the solution of \eqref{eq:rescaled} satisfies
  \begin{equation*}
    %\label{eq:estimH}
    \forall t>0,\qquad \fff(f(t))\leq \fff(f_0)e^{-\kappa t}
  \end{equation*}
\end{prop}

\preuve
We write (omitting the variable $t$ in the computations)
\begin{equation*}
\begin{split}
\ddt \fff(f(t)) & = \ddt \norm{f}^2 + \eps \ddt \seq{\Lambda_x^{-2} \D_x r, m}.
\end{split}
\end{equation*}
For the first term we notice that
\begin{equation} \label{premier}
\ddt \norm{f}^2 = \seq{L f, f} \leq -\norm{ h}^2
\end{equation}
 from the spectral gap property for $L$ in \eqref{zeinab}. For the second term we use the macroscopic equations and we get
 \begin{equation*}
\begin{split}
 \ddt \seq{\Lambda_x^{-2} \D_x r, m}
 & = \seq{\Lambda_x^{-2} \D_x r, \ddt m} + \seq{\Lambda_x^{-2} \D_x \ddt r, m}   \\
 & = -\seq{\Lambda_x^{-2}\D_x r, \D_x r}  + \seq{\Lambda_x^{-2} \D_x r, \Op_1(h)} + \seq{\Lambda_x^{-2} \D_x \Op_1(h), m}\\
 & \leq - \norm{\Lambda_x^{-1} \D_x r}^2 + C \norm{\Lambda_x^{-1} \Op_1(h)}\sep{ \norm{\Lambda_x^{-1} \D_x r)} + \norm{\Lambda_x^{-1} \D_x m}}
\end{split}
\end{equation*}
Now we use that $\norm{m} \leq \norm{h}$,  the preceding lemma and the Cauchy-Schwartz inequality  so that
\begin{equation*}
\begin{split}
 \ddt \seq{\Lambda_x^{-2} \D_x r, m}
 & \leq - \frac{1}{2} \norm{\Lambda_x^{-1} \D_x r}^2 + C \norm{h}^2.
\end{split}
\end{equation*}
The $L^2(\dmu)$ Poincar\'e inequality can be applied since $\seq{r} = \seq{f} = \seq{f_0} =  0$ and we get
\begin{equation} \label{second}
\begin{split}
 \ddt \seq{\Lambda_x^{-2} \D_x r, m}
 & \leq - \frac{1}{2} \frac{c_p}{c_p+1} \norm{r}^2 + C \norm{h}^2
\end{split}
\end{equation}
%with $\kappa_0 = \frac{c_p}{2(c_p + 1)}$.
Putting results from \eqref{premier} and \eqref{second} together gives
$$
\ddt \fff(f(t)) \leq - \norm{h}^2 - \frac{\eps}{2} \frac{c_p}{c_p+1} \norm{r}^2 + C \eps \norm{h}^2
$$
Now just taking  $\eps$ such that $C\eps \leq 1/2$ gives
\begin{equation*}
\begin{split}
\ddt \fff(f(t)) \leq - \frac {1}{2} \norm{h}^2 - \frac{\eps}{2} \frac{c_p}{c_p+1} \norm{r}^2 \leq - \frac{\eps}{2} \frac{c_p}{c_p+1} \norm{f}^2 \leq - \frac{\eps}{4} \frac{c_p}{c_p+1} \fff(f(t)).
\end{split}
\end{equation*}
This gives the result with $2\kappa = \frac{\eps}{4} \frac{c_p}{c_p+1}$. \fin

We can now easily prove the main result.

\preuve[ of Theorem \ref{hypocoercivite}]
We just have to notice that according to Lemma \ref{lem:equiv2} and proposition \ref{prop:decrexp2} we have
$$
\norm{f(t)}^2 \leq 2 \fff(f(t) \leq 2 \exp^{-2\kappa t} \fff(f_0) \leq 4 \exp^{-2\kappa t} \norm{f_0}^2.
$$
The proof is complete. \fin

%%%%%%%%%%%%%%%%%%%%%%%%%%%%%%%%%%%%%%%%%%%%%%%%%%%%%%%%%%%%%
\section{Short time regularization}

In this section we prove Theorem \ref{hypoellipticite} concerning diffusive kinetic equations.  The result of Theorem \ref{hypoellipticite} may seem to be very different from the one of Corollary \ref{hun}, since the last concerns exponential time decay, but the proof is in fact very similar. It is again base on the fundamental equality $\adf{\D_v, v\D_x} = \D_x$ and in fact we shall use computations done there in the proof here.

\subsection{A new entropy for  (FP) and the regularization effect}
Again we shall focus on a very simple case. We consider the  (FP) case with $d=1$, $V=0$ and $x\in \T$. The main step is to build a good entropy again:  for a function $f \in H^1(\dmu\dnu)$, we denote for C, D, E to be chosen later
$$
\ggg(t,f)   = C\norm{f}^2+Dt\norm{\partial_vf}^2+Et^2\seq{\partial_v f,\partial_x f}+t^3 \norm{\partial_xf}^2.
$$
 The aim of the following proposition is to show that $\ggg(t,f(t))$ is indeed a good entropy functional when $f$ is the solution of  (FP). In the following we use (without proof) that the equation is well posed in $H^1(\dmu\dnu)$ as we already did in the proof of Proposition \ref{prop:decrexp1}.

\begin{prop}\label{prop:decrexp3}
   There exits a constants $C$, $D$, $E$  such that for all $f_0 \in H^1(\dmu\dnu)$
  solution of \eqref{eq:rescaled} the modified entropy $\ggg(f(t))$  satisfies
  \begin{equation*}
    %\label{eq:estimH}
    \forall t\in [0,1],\qquad \ggg(t,f(t))\leq \ggg(0,f_0)
  \end{equation*}
\end{prop}

\preuve
We just have to show that
$$
\ddt \ggg(t,f(t)) \leq 0.
$$
Several terms are involved and we write
\begin{equation*} %\label{decomp}
\begin{split}
\ddt \ggg(t,f(t)) = &   C \ddt \norm{f}^2+Dt \ddt \norm{\partial_vf}^2+Et^2 \ddt \seq{\partial_v f,\partial_x f}+t^3 \ddt \norm{\partial_xf}^2 \\
& +D\norm{\partial_vf}^2+2Et\seq{\partial_v f,\partial_x f}+3t^2 \norm{\partial_xf}^2.
\end{split}
\end{equation*}
We use freely the computations done in the proof of Proposition \ref{prop:decrexp1} and we get
\begin{equation*} %\label{decomp}
\begin{split}
\ddt \ggg(t,f(t)) \leq  &   -2C \norm{\D_v f}^2  -2 Dt \norm{(-\D_v+v) \D_v f}^2 - E t^2 \norm{\D_x f}^2 - 2t^3 \norm{\D_x \D_v f}^2 \\
& -2 D t \seq{\D_x f, \D_v f} +2 E t^2 \seq{(-\D_v + v) \D_v f, \D_x \D_v f} - E t^2 \seq{\D_x f, \D_v f} \\
& +D\norm{\D_vf}^2+2Et\seq{\D_v f,\D_x f}+3t^2 \norm{\D_xf}^2,
\end{split}
\end{equation*}
where we put on the first line all terms with a sign, on the second line all term without sign coming from the derivative in time, and in the last line the terms for which the derivation was done only on the prefactor in powers of $t$.

Now we use that we are on $t \in [0,1]$ and we get
\begin{equation*} %\label{decomp}
\begin{split}
\ddt \ggg(t,f(t)) \leq  &   -2C \norm{\D_v f}^2  -2 Dt \norm{(-\D_v+v) \D_v f}^2 - E t^2 \norm{\D_x f}^2 - 2t^3 \norm{\D_x \D_v f}^2 \\
& +( 2D + 3E) t |\seq{\D_x f, \D_v f}| +2 E t^2 |\seq{(-\D_v + v) \D_v f, \D_x \D_v f}|  \\
& +D\norm{\D_vf}^2+3t^2 \norm{\D_xf}^2.
\end{split}
\end{equation*}
We use Cauchy-Schwartz inequality in the middle line and we get
\begin{equation*}
\left\{
\begin{array}{l}
( 2D + 3E) t |\seq{\D_x f, \D_v f}| \leq ( 2D + 3E)^2 \norm{\D_v f}^2 + t^2 \norm{\D_x f}^2 \\
2 E t^2 |\seq{(-\D_v + v) \D_v f, \D_x \D_v f}|  \leq  E^2 t \norm{(-\D_v+v) \D_v f}^2 + \norm{\D_x \D_v f}^2
\end{array}
\right.
\end{equation*}
This provides the following result
\begin{equation*} %\label{decomp}
\begin{split}
\ddt \ggg(t,f(t)) \leq  &   (-2C + (2D + 3E)^2 +D) \norm{\D_v f}^2 + (-2 D + E^2) t \norm{(-\D_v+v) \D_v f}^2 \\
 & + (-E +3 +1)t^2 \norm{\D_x f}^2 + (- 2 +1) t^3 \norm{\D_x \D_v f}^2
\end{split}
\end{equation*}
Now $(-2+1) \leq 0$, then we choose $E$ such that $(-E+4) \leq 0$, next $D$ such that $(-2D+ E^2) \leq 0$ and at the end $C$ such that $(-2C + (2D + 3E)^2 +D) \leq 0$ and we get the result.

\preuve[of theorem \ref{hypoellipticite}]
The proof is direct. Notice first that $\ggg(0,f(0)) = C \norm{f_0}^2$. From the preceding Proposition we get that for all $t \geq 0$ and $f_0 \in H^1(\dmu\dnu)$, we have
$$
Dt\norm{\D_v \exp^{-tP} f_0}^2 \leq C \norm{f_0}^2 \qquad \textrm { and } \qquad t^3\norm{\D_x \exp^{-tP} f_0}^2 \leq C \norm{f_0}^2.
$$
Using the density of $H^1(\dmu\dnu)$  in $L^2(\dmu\dnu)$,  we get for all $t >0$ operators $\D_v e^{-tP}$ and $\D_x \exp^{-tP} $ are $L^2$ bounded with respective bound $\sqrt{\frac{C}{D}}t^{-1/2}$ and $\sqrt{C} t^{-3/2}$. The proof is complete. \fin

\subsection{The fractional Kolmogorov equation}

This section is devoted to the proof of Theorem \ref{FKEthm} concerning the fractional Kolmogorov equation (FK).   We follow the lines of the proof given in \cite{Her07} and \cite{HTT16a}. Since we will work below on the Fourier side, we work from now on  in the (complex) $L^2$. We denote again $\norm{\cdot}$ is the usual $L^2$ norm and $\seq{ \cdot,\cdot}$ the usual (complex) $L^2$ scalar product. Note anyway that  we deal with a real operator, and that any solution of the evolution equation \eqref{FKE} is real if the corresponding initial datum is. For a function $f \in H^s$ we introduce a new  adapted entropy functional $\kkk$ depending on time and defined  for all $t\geq 0$ by
$$
\kkk(t,f) = C\norm{f}^2 + D t \norm{\Lambda_v^{s-1} \D_v f}^2 + E t^{1+s} \Re \seq{\Lambda_v^{s-1} \D_v f, \Lambda_x^{s-1} \D_x f} + t^{1+2s} \norm{ \Lambda_x^{s-1} \D_x f}^2
$$
for large constants $C$, $D$, $E$ to be chosen later. The first step in the study is to show that
that $\kkk$ is indeed non-negative. The lemma below shows in addition that for all $t > 0$,
$\kkk(t,f)$ controls the $H^s$ norm.

\begin{lem} \label{hh1}
If $E^2 \leq D$ then for all  $t\geq 0$ and $f \in H^{s}$ we have $\kkk(t,f) \geq 0$. Precisely we have
\begin{equation*}
0 \leq  C\norm{f}^2 + \frac{D}{2} t \norm{\Lambda_v^{s-1} \D_v f}^2  + \frac{1}{2} t^{1+2s} \norm{ \Lambda_x^{s-1} \D_x f}^2
\leq \kkk(t,f)% \leq  C\norm{f}^2 + \frac{3}{2} D t \norm{\Lambda_v^{s-1} \D_v f}^2 + \frac{3}{2} t^{1+2s} \norm{ \Lambda_x^{s-1} \D_x f}^2
\end{equation*}
\end{lem}
\preuve The proof is direct using the time-dependant Cauchy-Schwartz inequality
$$
E t^{s} \abs{ \sep{\Lambda_v^{s-1} \D_v f, \Lambda_x^{s-1} \D_x f}}  \leq
\frac{E^2}{2}  \norm{\Lambda_v^{s-1} \D_v f}^2  + \frac{1}{2} t^{2s} \norm{ \Lambda_x^{s-1} \D_x f}^2.
$$
\fin

As before, the main ingredient in the proof of Theorem \ref{FKEthm} is the following commutation equality
$$
\adf{ \D_v, v \D_x} = \D_x.
$$
In the same spirit we shall need later the following lemma
giving formulas for  slightly modified commutators. We denote from now on $X_0 = v\D_x$ the Vlasov operator, so that the preceding fundamental equality reads $ \adf{ \D_v, X_0} = \D_x$.

\begin{lem} \label{commm}
We have $$ \adf{ \Lambda_v^{s-1} \D_v, X_0} = { (1-s\Delta_v)}\Lambda_v^{s-3}  \D_x
  $$
  and % \textrm{   and }
$$
\adf{ \Lambda_v^{s-1} \D_v, \Lambda_v^{2s}} = \adf{ \Lambda_x^{s-1} \D_x, \Lambda_v^{2s}} = \adf{ \Lambda_x^{s-1} \D_x, X_0} = 0.$$
\end{lem}

\preuve For the three last equalities, the result is immediate since differentiation in velocity and spatial direction commute. Let us deal with the first one. We check that
that the commutator $\adf{ \Lambda_v^{s-1} \D_v, X_0}$ is in fact a Fourier multiplier whose symbol reads
\begin{equation*} %\label{calcc}
\begin{split}
\sigma \sep{ \adf{ \Lambda_v^{s-1} \D_v, X_0}} & = \frac{1}{i} \set{ (1+\eta^2)^{\frac{s-1}{2}} i\eta, iv \xi}
\end{split}
\end{equation*}
where we denote by $\set{ \cdot, \cdot}$ the Poisson bracket of two functions. In the preceding Fourier formalism, we used that
\begin{eqnarray*}
\sigma( \D_v) = i\eta, \quad \sigma(\D_x) = i\xi, \qquad \sigma(X_0) = iv\xi \\
\sigma(\Lambda_v^\alpha) = \seta^\alpha, \quad \sigma( \Lambda_v^\beta) = \seta^\beta, \quad \sigma( -\Delta_v) = \eta^2.
\end{eqnarray*}
We have then
\begin{multline}
\sigma \sep{ \adf{ \Lambda_v^{s-1} \D_v, X_0}}  = i \xi \set{ (1+\eta^2)^{\frac{s-1}{2}} \eta, v} \\
 =  i\xi \sep{ (s-1) \eta^2 \seta^{s-3} + \seta^{s-1}}
 =  (1+s\eta^2)\seta^{s-3} i\xi.
\end{multline}
Coming back in the non-Fourier side, we get the result.
\fin

We now show that $\kkk$ is indeed a Lyapunov function. Let $f$ be a solution of
\begin{equation} \label{FKEbis}
 \D_{t}f+v\D_{x}f + \Lambda_v^{2s}f =0
 \end{equation}
 with (real) initial data $f_0$. We have then
\begin{prop} \label{deriv}  For well chosen constants $C$, $D$ and $E$ we have
$$
\ddt \kkk(t,f(t)) \leq 0.
$$
\end{prop}

\preuve Using the previous lemma, we shall compute the time derivative of each of the terms appearing in the definition of $\kkk$. For convenience  we introduce the operator associated the the Kolmogorov equation
$$
P = X_0 + \Lambda_v^{2s}.
$$
so that $f$ satisfies $\D_t f+ Pf = 0 $.  We first notice that
 \begin{equation*} %\label{ordre0}
 \ddt \norm{f}^2 = -2\Re \seq{Pf,f} = -2\Re \seq{ (X_0 + \Lambda^{2s}_v)f,f} =  -2  \seq{\Lambda_v^{2s} f,f}
 \end{equation*}
 since $X_0$ is skewadjoint. Using Parseval formula on the right side we get that the first term in the derivative of $\kkk$ is
  \begin{equation} \label{f}
 \ddt C\norm{f}^2 =  -\big\langle \underbrace{2 C \seta^{2s}}_{I} \hf, \hf \big\rangle.
 \end{equation}
 Note that this term is negative.
For the second term in the derivative of $\kkk$ we have
\begin{equation*}
\ddt \sep{  t\norm{\Lambda_v^{s-1} \D_v f}^2 } =  \norm{\Lambda^{s-1}_v \D_v f}^2 + t \ddt \seq{\Lambda^{s-1}_v \D_v f , \Lambda^{s-1}_v \D_v f}.
 \end{equation*}
The derivative in the last term writes
\begin{equation*}\begin{split}
\ddt \norm{\Lambda^{s-1}_v \D_v f}^2  = & -2 \Re \seq{\Lambda_v^{s-1} \D_v P f , \Lambda_v^{s-1} \D_v f}  \\
  = &  -2 \Re \seq{ \Lambda_v^{s-1} \D_v \Lambda_v^{2s} f , \Lambda_v^{s-1} \D_v f} -2 \Re \seq{\Lambda_v^{s-1} \D_v  X_0 f , \Lambda_v^{s-1} \D_v f} \\
   = & -2 \Re \seq{ \Lambda_v^{2s} \Lambda_v^{s-1} \D_v  f , \Lambda_v^{s-1} \D_v f} -2 \Re \seq{ X_0 \Lambda_v^{s-1} \D_v   f , \Lambda_v^{s-1} \D_v f} \\
  & -2 \Re \seq{\adf{ \Lambda_v^{s-1} \D_v , \Lambda_v^{2s}} f , \Lambda_v^{s-1} \D_v f} -2 \Re \seq{ \adf{ \Lambda_v^{s-1} \D_v, X_0  }  f , \Lambda_v^{s-1} \D_v f} \\
  = & -2 \Re \seq{\Lambda_v^{2s} \Lambda_v^{s-1} \D_v  f , \Lambda_v^{s-1} \D_v f}  - 2 \Re \seq{  (1-s\Delta_v)\Lambda_v^{s-3} \D_x f , \Lambda_v^{s-1} \D_v f} \\
 \end{split}
 \end{equation*}
 where  we used that $X_0$ is skewadjoint and the commutation expressions in  Lemma \ref{commm}. Writing the right-hand side on the Fourier side gives
 \begin{equation*}\begin{split}
\ddt \norm{\Lambda^{s-1}_v \D_v f}^2
  & =   2 \seq{ \Lambda_v^{4s-2} \Delta_v  f ,  f}  + 2 \Re \seq{ (1-s\Delta_v) \Lambda_v^{2s-4} \D_x \D_v f , f} \\
  & = - 2 \seq{ \seta^{4s-2} \eta^2 \hf, \hf} -2 \seq{ (1+s\eta^2)\seta^{2s-4} \eta \xi \hf, \hf} \\
  & \leq - 2 \seq{ \seta^{4s}  \hf, \hf} +2 \seq{ \seta^{4s-2}  \hf, \hf} + 2 \seq{ \seta^{2s-1}  \sxi \hf, \hf}
 \end{split}
 \end{equation*}
 since $\abs{  (1+s\eta^2)\seta^{2s-4} \eta \xi} \leq \seta^{2s-1}  \sxi$.
 The second term in $\kkk$ therefore satisfies
\begin{multline} \label{vvf}
 \ddt \sep{ D t \norm{\Lambda^{s-1}_v \D_v f}^2 } \\ \leq \big\langle \big( \underbrace{D \seta^{2s}}_{i} - \underbrace{2Dt\seta^{4s}}_{II} + \underbrace{2Dt\seta^{4s-2}}_{ii} + \underbrace{2 Dt  \seta^{2s-1} \sxi}_{iii}\big) \hf, \hf \big\rangle.
 \end{multline}
 We note that the term corresponding to II is negative, and that the three other ones are positive.
 Now we deal with the third term in the derivative of $\kkk$:
\begin{multline*}
\ddt \sep{  t^{1+s} \Re\seq{\Lambda_v^{s-1} \D_v f, \Lambda_x^{s-1} \D_x f} } \\
= (1+s) t^{s} \Re \seq{\Lambda_v^{s-1} \D_v f, \Lambda_x^{s-1} \D_x f} + t^{1+s} \ddt  \Re \seq{\Lambda_v^{s-1} \D_v f, \Lambda_x^{s-1} \D_x f} .
 \end{multline*}
The derivative in the last term writes
\begin{equation*}\begin{split}
\ddt \Re\seq{\Lambda_v^{s-1} \D_v f , \Lambda_x^{s-1} \D_x f}  = & - \Re \seq{\Lambda_v^{s-1} \D_v P f , \Lambda_x^{s-1} \D_x f}- \Re \seq{\Lambda_v^{s-1} \D_v  f , \Lambda_x^{s-1} \D_x P f}  \\
  = &  -2 \Re \seq{\Lambda_v^{2s} \Lambda_v^{s-1} \D_v  f , \Lambda_x^{s-1} \D_x f} \\
    & - \Re \seq{\Lambda_v^{s-1} \D_v  X_0 f , \Lambda_x^{s-1} \D_x f}- \Re \seq{\Lambda_v^{s-1} \D_v  f , \Lambda_x^{s-1} \D_x X_0 f} \\
    = &  -2 \Re \seq{ \Lambda_v^{2s} \Lambda_v^{s-1} \D_v  f , \Lambda_x^{s-1} \D_x f} \\
    & - \Re \seq{\adf{ \Lambda_v^{s-1} \D_v , X_0} f , \Lambda_x^{s-1} \D_x f}- \Re \seq{\Lambda_v^{s-1} \D_v  f , \adf{ \Lambda_x^{s-1} \D_x, X_0 } f} \\
    & - \Re \seq{ X_0 \Lambda_v^{s-1} \D_v   f , \Lambda_x^{s-1} \D_x f}- \Re \seq{ \Lambda_v^{s-1} \D_v  f , X_0 \Lambda_x^{s-1} \D_x  f}.
    \end{split}
    \end{equation*}
    Now we use again that $X_0$ is skewadjoint and observe that it implies that the sum of the last two terms is zero by compensation. The previous term is also zero since the commutator inside is zero. With Lemma \ref{commm} we obtain
    \begin{multline*}
    \ddt \Re \seq{ \Lambda_v^{s-1} \D_v f , \Lambda_x^{s-1} \D_x f} \\
    =  -2 \Re \seq{ \Lambda_v^{2s} \Lambda_v^{s-1} \D_v  f , \Lambda_x^{s-1} \D_x f}
    - \Re \seq{ (1-s\Delta_v) \Lambda_v^{s-3} \D_x f , \Lambda_x^{s-1} \D_x f}.
 \end{multline*}
  Writing the right-hand side on the Fourier side gives then
 \begin{equation*}\begin{split}
\ddt \Re \seq{ \Lambda_v^{s-1} \D_v f , \Lambda_x^{s-1} \D_x f}
  & = -  \seq{ \seta^{3s-1} \sxi^{s-1} \eta \xi \hf, \hf}  -  \seq{ (1+s\eta^2)\seta^{s-3} \sxi^{s-1}  \xi^2 \hf, \hf} \\
   & \leq  \seq{ \seta^{3s} \sxi^{s}   \hf, \hf }
  -s \seq{ \seta^{s-1} \sxi^{s-1} \xi^2   \hf, \hf } \\
    & \leq  \seq{ \seta^{3s} \sxi^{s}   \hf, \hf}
  -s \seq{ \seta^{s-1} \sxi^{s+1}   \hf, \hf}
  + s \seq{ \seta^{s-1} \sxi^{s-1}   \hf, \hf}.
 \end{split}
 \end{equation*}
 %Now we use that
% $$
% \seta^{s-1} \sxi^{s-1} (\eta^2+\alpha) \xi^2 \leq C_s \seta^{s-1} \sxi^{s+1} + \seta^{s-1} \sxi^{s+1}
% $$
 We therefore get that the third term in $\kkk$ therefore satisfies
\begin{multline} \label{vxf}
 \ddt \sep{ E  t^{1+s} \Re\sep{\Lambda_v^{s-1} \D_v f, \Lambda_x^{s-1} \D_x f} } \\ \leq \big\langle \big( \underbrace{E(s+1)t^s \seta^{s} \sxi^{s} }_{iv}   + \underbrace{E  t^{1+s}\seta^{3s}\sxi^{s}}_{v} \\ - \underbrace{E s t^{1+s} \seta^{s-1}\sxi^{s+1}}_{III} + \underbrace{E s t^{1+s}  \seta^{s-1}  \sxi^{s-1}}_{vi}  \big) \hf, \hf \big\rangle.
 \end{multline}
 We note that the term corresponding to III is negative, and that the three other ones are positive.

We can now deal with the last term in the derivative of $\kkk$.
We write

\begin{equation*}
\ddt \sep{  t^{1+2s}\norm{\Lambda_x^{s-1} \D_x f}^2 } =  (1+2s) t^{2s} \norm{\Lambda_x^{s-1} \D_x f}^2 + t^{1+2s} \ddt \norm{\Lambda_x^{s-1} \D_x f }^2.
 \end{equation*}
The last time derivative  writes
\begin{equation*}\begin{split}
\ddt \norm{\Lambda_x^{s-1} \D_x f }^2  = & -2 \Re \seq{ \Lambda_x^{s-1} \D_x P f , \Lambda_x^{s-1} \D_x f}  \\
  = &  -2 \Re \seq{ \Lambda_v^{2s} \Lambda_x^{s-1} \D_x  f , \Lambda_x^{s-1} \D_x f} -2 \Re \seq{ X_0 \Lambda_x^{s-1} \D_x  f , \Lambda_x^{s-1} \D_x f} \\
  = & -2 \Re \seq{ \Lambda_v^{2s} \Lambda_x^{s-1} \D_x  f , \Lambda_x^{s-1} \D_x f}.
 \end{split}
 \end{equation*}
 We used here the last commutations properties in Lemma \ref{commm} and again that $X_0$ is skewadjoint. Writing the right-hand side on the Fourier side gives
 \begin{equation*}\begin{split}
\ddt \norm{\Lambda_x^{s-1} \D_x f }^2
   = - 2\seq{ \seta^{2s} \Lambda_x^{2s-2} \xi^2 \hf, \hf}
  & = - 2\seq{\seta^{2s} \sxi^{2s}  \hf, \hf}) +2\seq{\seta^{2s} \sxi^{2s-2} \hf, \hf}
 \end{split}
 \end{equation*}
 The fourth term in the derivative of  $\kkk$ therefore satisfies
\begin{multline} \label{xxf}
 \ddt \sep{ t^{1+2s} \norm{\Lambda_x^{s-1} \D_x f }^2 }  \\
 \leq \big\langle \big( \underbrace{(1+2s) t^{2s} \sxi^{2s}}_{vii} - \underbrace{2t^{1+2s} \seta^{2s}\sxi^{2s}}_{IV} + \underbrace{2t^{1+2s} \seta^{2s}\sxi^{2s-2}}_{viii}\big) \hf, \hf \big\rangle.
 \end{multline}
 We note that the term corresponding to IV is negative, and that the other ones are positive.

\bigskip
Now we look at the different
terms appearing in formulas (\ref{f}-\ref{xxf}). We want to show that with a good choice of constants $C$, $D$ and $E$, the corresponding sum is negative, and therefore $\kkk$ is indeed a Lyapunov functional. We shall study each non-negative term (small letters $i$ to $viii$ ) and show that they can be controlled by combinations of term $I$ to $IV$, using essentially the H\"older inequality in $\R^2$. We restrict the study to $t \in [0,1]$.

The terms $(i)$  and $(ii)$ are immediately bounded by $I/10$ if
\begin{equation} \label{condi}
2D \leq 2C/10
\end{equation}
 since $s\leq 1$.
The term $(iii)$ is more involved. We check that for any $\eps_{iii} >0$
$$
t \seta^{2s-1} \sxi \leq \eps_{iii}^{-1} \seta^{2s} + \eps_{iii}^s t^{1+s} \seta^{s-1}\sxi^{s+1}.
$$
Multiplying this inequality by $2D$ implies that $(iii) \leq I/10 + III/10$ if the following conditions are satisfied
\begin{equation} \label{condiii}
\eps_{iii}^{-1} 2D \leq 2C/10,  \qquad \eps_{iii}^s 2D \leq Es/10.
\end{equation}

Now we deal with the term $(iv)$. We first check that for any $\eps_{iv} >0$
$$
t^s \seta^{s} \sxi^s \leq \eps_{iv}^{-1} \seta^{2s} + \eps_{iv}^{1/s} t^{1+s} \seta^{s-1}\sxi^{s+1}.
$$
Multiplying this inequality by $E(s+1)$ implies that $(iv) \leq I/10 + III/10$ if the following conditions are satisfied
\begin{equation} \label{condiv}
\eps_{iv}^{-1} E(s+1) \leq 2C/10,  \qquad \eps_{iv}^{1/s} E(s+1) \leq Es/10.
\end{equation}
For the term $(v)$ we also have to give a refined estimate.
We first check that for any $\eps_{v} >0$
$$
t^{1+s} \seta^{3s} \sxi^s \leq \eps_{v}^{-1} t \seta^{4s} + \eps_{v} t^{1+2s} \seta^{2s}\sxi^{2s}.
$$
Multiplying this inequality by $E$ implies that $(v) \leq II/10 + IV/10$ if the following conditions are satisfied
\begin{equation} \label{condv}
\eps_{v}^{-1} E \leq 2D/10,  \qquad \eps_{v} E \leq 2/10.
\end{equation}
The term $(vi)$ is easily handled since $s\leq 1$, and we directly get that
$(vi) \leq I/10$ if
\begin{equation} \label{condvi}
Es \leq 2C/10.
\end{equation}
Now we study the term $(vii)$. We first notice that for any $\eps_{vii} >0$
$$
t^{2s}  \sxi^{2s} \leq \eps_{vii}^{-1}  \seta^{2s} + \eps_{vii}^{\frac{1-s}{2s}} t^{1+s} \seta^{s-1}\sxi^{s+1}.
$$
Multiplying this inequality by $(1+2s)$ implies that $(vii) \leq I/10 + III/10$ if the following conditions are satisfied
\begin{equation} \label{condvii}
\eps_{vii}^{-1} (1+2s) \leq 2C/10,  \qquad \eps_{vii}^{\frac{1-s}{2s}} (1+2s) \leq Es/10.
\end{equation}
To finish The term $(viii)$ is also easily handled since $s\leq 1$, and we directly get that
$(viii) \leq I/10$ if
\begin{equation} \label{condviii}
2 \leq 2C/10.
\end{equation}

Now we can do the synthesis and check that we can choose (in order  of reverse appearance) the constants $C$, $D$, $E$ and the small constants
$\eps_{iii}$, $\eps_{iv}$, $\eps_{v}$ and $\eps_{vii}$ such that
conditions (\ref{condi}-\ref{condviii}) are satisfied. Note that $D$ and after that $C$  can be taken arbitrarily larger at the end of this procedure. We obtain therefore that
\begin{equation} \label{synthpartielle}
\ddt \kkk(t, f(t)) \leq -\frac{1}{10} \seq{(I + II + III + IV)\hf, \hf} \leq 0
\end{equation}
and the proof is complete. \fin

%In fact the proof before (see \eqref{synthpartielle}) gives the much stronger following result of fundamental use.
%\begin{lem} For well chosen (arbitrarly large) constants $C$, $D$ and $E$ there exists $\kappa >0$ such that for all $t\geq 0$
%$$
%\ddt \hhh(t) \leq -\kappa \hhh(t).
%$$
%\end{lem}
%
%\preuve Using \eqref{synthpartielle} and Lemma \ref{hh1} it is sufficient
%to prove that for
%\begin{equation} \label{IIIIII}
%2\kappa \sep{ C\norm{f}^2 + \frac{3}{2} D t \norm{\Lambda_v^{s-1} \D_v f}^2 + \frac{3}{2} t^{1+2s} \norm{ \Lambda_x^{s-1} \D_x f}^2 } \leq  \frac{1}{10} \sep{  (I + II + III + IV) \hf, \hf}
%\end{equation}
%Recall that all the terms $I$ to $IV$ are positive. We compare each term on the Fourier side. Indeed \eqref{IIIIII} is clearly implied by the following inequality  :
%$$
%2\kappa C + 3 \kappa D t \seta^{2s-2} \eta^2 + 3 \kappa t^{1+2s}\sxi^{2s-2} \xi^2 \leq \frac{1}{10} \sep{ 2C \seta^{2s} + Dt \seta^{4s} + Est^{1+2s} \seta^{s-1} \sxi{s+1} + t^{1+2s} \seta^{2s} \sxi{2s}}
%$$

Now we are able to conclude the proof of the main result Theorem~\ref{FKEthm} concerning  the fractional Kolmogorov equation.

\preuve[of Theorem \ref{FKEthm}]
We prove first the result for $r=0$.
Let $C$, $D$ and $E$ be constants given by Propositions \ref{deriv} and  \ref{hh1} and let us take $f_0 \in \sss$.
From Lemma \ref{deriv} we first get that for all $t \in [0,1]$
$$
\kkk(t,f(t)) \leq \kkk(0,f_0) = C \norm{f_0}^2.
$$
Now using Lemma \ref{hh1}, we get in particular
$$
\frac{D}{2} t\norm{ \Lambda_v^s f}^2 \leq C\norm{f}^2 + \frac{D}{2} t \norm{\Lambda_v^{s-1} \D_v f}^2 \leq \kkk(t,f(t)) \leq C\norm{f_0}^2
$$
and this implies the result for the velocity regularization.
Similarly using again \ref{hh1}, we have
$$
\frac{1}{2} t^{1+2s} \norm{ \Lambda_x^{s}  f}^2 \leq C\norm{f}^2 +\frac{1}{2}t^{1+2s} \norm{ \Lambda_x^{s-1} \D_x f}^2 \leq \kkk(t,f(t)) \leq C\norm{f_0}^2
$$
and this gives the regularization result for $r=0$ in the spatial direction.

Now for $r\in \R$ we just use the fact that $P$ commutes with
$\Lambda_x^r$ which implies that for $f$  solution of $\D_t f + Pf = 0$ with initial data $f_0$,   $\Lambda_x^r f$ is the solution of
$\D_t \Lambda_x^r f + P \Lambda_x^r f = 0$ with initial data $\Lambda_x^r f_0$. We can therefore apply the result on $L^2$ to $\Lambda_x^r f$ and this directly that
$$
\frac{D}{2} t\norm{ \Lambda_v^s\Lambda_x^r f}^2  \leq C\norm{\Lambda_x^r f_0}^2 \textrm{ and } \frac{1}{2} t^{1+2s} \norm{ \Lambda_x^{s+r}  f}^2 \leq  C\norm{\Lambda_x^r f_0}^2.
$$
This gives the estimates for any $r\in \R$ and $f_0 \in \sss$. The general result for initial data in the corresponding spaces follows by density of $\sss$. The proof is complete. \fin

\section{Application to a mollified Vlasov-Poisson-Fokker-Planck system}

In this  subsection we apply the results of the previous two sections to the  simple non-linear system (MVPFP) in a perturbative situation and prove Theorem \ref{MVPFPbis}.
We study here the $1d$-model where $V=0$, $(x,v)\in\T\times\R$. Recall that it reads
\begin{equation*}
\left\{
\begin{array}{l}
 \D_t F + v \D_x F +\Enl(F) \D_v (F-\mmm) = \D_v(\D_v+v)F,\\
 \pm \Enl(F) = K * (\rho-1), \qquad \textrm{ with }  \rho = \int F \dv \\
 F|_{t=0} = F_0, \qquad \iint F_0 \dx\dv = 1,
 \end{array}
 \right.
 \end{equation*}
 where $K$ is supposed to be smooth and in $ L^\infty(\dx)$.
  We can directly check that the global Maxwellian $\mmm = \nu$ is a stationary solution and the question of the existence  and the trend to the equilibrium in a perturbative context raises naturally. We shall below prove both simultaneously.

  \bigskip
 For doing this, we shall again change the variable and pose
 $$
 F = \mmm + \mmm f.
 $$
 (Recall that here $\mmm(x,v) = \nu(v)$). Integrating the equation gives that the mass is conserved which reads in the new variables $\seq{f(t)} = \seq{f_0} = 0$. The new equation reads then
 \begin{equation*}
\left\{
\begin{array}{l}
 \D_t f + v \D_x f - E(f)(-\D_v+v)f +(-\D_v+v) \D_v f = 0 \\
 \pm E(F) = K * r, \qquad \textrm{ with }  r(x) = \int f \dnu \\
 f|_{t=0} = f_0, \qquad \seq{f_0} = 0
 \end{array}
 \right.
 \end{equation*}
The non-linear term is then considered as a perturbation and we write
$$
 \D_t f + v \D_x f +(-\D_v+v) \D_v f =  E(f)(-\D_v+v)f
 $$
 One of the main difficulty is that although there is a spectral gap for the operator (FP) appearing in the left, it cannot absorbed the term on the right which is of order $1$ in velocity. For bypassing this problem we shall use the regularizing properties. Writing $Pf = v \D_x f +(-\D_v+v) \D_v f $ and using Duhamel formula we get
 $$
 f(t) = \exp^{-tP} f_0 + \int_0^t \exp^{-(t-s)P } E(f(s)) (-\D_v + v) f(s)  \ds.
 $$
 (we omit the dependance on $(x,v)$ variables).
 Remembering that $E(f)$ only depends on the variable $(t,x)$ and using this in the  Duhamel formula we get
 \begin{equation} \label{duhamel}
 \left\{
 \begin{array}{l}
 f(t) = \exp^{-tP} f_0 + \int_0^t E(f(s)) \exp^{-(t-s)P}  (-\D_v + v) f(s)  \ds, \\
 \pm E(f) =  K * \int f \dnu, \qquad f|_{t=0} = f_0, \qquad \seq{f_0} = 0
 \end{array}
 \right.
 \end{equation}
 Any couple $(f, E)$ satisfying \eqref{duhamel} will be called a \it mild \rm solution of the mollified Vlasov-Poisson-Fokker-Planck system.

 We first prove a fundamental Lemma in the spirit of the Sections 2 and 3 about the operator appearing in the Duhamel term:
 \begin{lem} \label{together}
 There exists a constant such that For all $g \in L^2(\dmu\dnu)$
  $$
 \norm{\exp^{-tP}(-\D_v + v)g} \leq C (t^{-1/2} + 1) \exp^{-t\kappa} \norm{g}.
 $$
 \end{lem}

 \preuve We use the results of the two previous sections.
 From Theorem \ref{hypoellipticite} we know that for $g \in L^2(\dmu\dnu)$ we have for all $t\in ]0,1]$
 $$
 \norm{\D_v \exp^{-tP}g} \leq C t^{-1/2}\norm{g}.
 $$
 Now all the computation before could be straight forwardly adapted to the adjoint $P^* = -X_0-L$ of $P$ and we also have
 $$
 \norm{\D_v \exp^{-tP^*}g} \leq C t^{-1/2}\norm{g}.
 $$
 Using $\D_v^* = (-\D_v+v)$ and taking the adjoint operator in the previous formula yields
 \begin{equation} \label{zeinab2}
 \norm{\exp^{-tP}(-\D_v + v)g} \leq C t^{-1/2}\norm{g}.
 \end{equation}
 Now  for all $g \in H^1(\dmu\dnu)$, we have $\seq{(-\D_v + v)g} = 0$ since
 $$
 \seq{(-\D_v + v)g} = \seq{(-\D_v + v)g, 1} = \seq{g, \D_v 1} = 0.
 $$
For all $t\geq 1$ we therefore have
 \begin{multline}
 \norm{\exp^{-tP}(-\D_v + v)g} =\norm{\exp^{-(t-1)P -P}(-\D_v + v)g} \\ \leq C \exp^{-(t-1)\kappa} \norm{\exp^{-P}(-\D_v + v)g} \leq C' \exp^{-t\kappa} \norm{g_0},
 \end{multline}
 where we used the  regularization result stated in \eqref{zeinab2} for the last inequality.
 The proof is complete. \fin

%By solution of the non linear equation we mean mild solution, that is $f$ is a solution if
%$f$ satisfies the Duhamel formula

We also need a second Lemma to control the field.

\begin{lem} \label{estimE} For all $\phi \in L^2(\dmu\dnu)$ we have $\norm{E(\phi)}_{L^\infty(\dx)} \leq \norm{K}_{L^\infty(\dx)} \norm{\phi}_{L^2(\dmu\dnu)}$.
\end{lem}

\preuve We just write that
\begin{equation*}
\begin{split}
\norm{E(\phi)}_{L^\infty} & = \norm{ K * \int \phi \dnu}_{L^\infty} \leq \norm{K}_{L^\infty} \norm{\int \phi \dnu}_{L^1(d\mu)} \\
& \leq \norm{K}_{L^\infty} \norm{\phi}_{L^1(\dmu\dnu)} \leq \norm{K}_{L^\infty} \norm{\phi}_{L^2(\dmu\dnu)}
\end{split}
\end{equation*}
where the last inequality come form the fact that $\dmu\dnu$ is a probability measure.
\fin

% Using this in \eqref{duhamel} and taking the norm we get
% $$
% \norm{f(t)} \leq \exp^{-\kappa t} \norm{f_0} + \int_0^t \max((t-s)^{-1/2}, 1) e^{-\kappa (t-s)} \norm{f(s)} \norm{E(f(s))}_{L^\infty}.
% $$
 %The following estimate comes easily
% $$
% \norm{E(f(s))}_{L^\infty} \leq \norm{k}_{L^\infty}\norm{r}_{L^1(\dmu)} \leq \sqrt{2\pi} \norm{k}_{L^\infty}\norm{r}_{L^2(\dmu)}

Theorem \ref{MVPFPbis} is a direct consequence of the  following Theorem, that we shall prove with the use of a Picard scheme and the two preceding Lemmas.

 \begin{thm} \label{MVPFP} Suppose $d=1$, $V=0$ and work on $\T$ in the space variable. Then there exists $\eps_0$, $C_0$ and $\kappa_0$ explicitly computable such that the following happens. Suppose $\seq{f_0} = 0$.
 \begin{enumerate}
 \item If $\norm{f_0} \leq \eps_0$ then there exists a unique mild solution $(f,E) \in \ccc^0(\R_+, L^2(\dmu\dnu) ) \times \ccc^0(\R^+, L^\infty(\dx))$ to the mollified Vlasov-Poisson-Fokker-Planck system \eqref{duhamel}.
     \item We have in addition for all $t\geq 0$,
     $$
     \norm{f(t)}_{L^2(\dmu\dnu)} \leq  C_0 \exp^{-\kappa_0 t} \quad \textrm{ and } \quad \norm{E(f(t))}_{L^\infty(\dx)} \leq  C_0 \exp^{-\kappa_0 t}
     $$
     \end{enumerate}
     \end{thm}

\preuve  Let us consider $f_0$ such that $\norm{f_0} \leq \eps_0$ to be fixed later. We  pose (omitting the space and velocity variables)
$$
f(t) = \exp^{-tP} f_0 + g(t), \qquad G(t) = K * \int g \dnu \qquad \Gl(t) = K * \int \exp^{-tP} f_0 \dnu
$$
We build a solution with a fixed point argument for the following map
$\Phi = (\Phi_1, \Phi_2)$ given by
\begin{equation*}
\Phi_{1}(g,G)(t)= \int_0^t \big(\Gl(s)+G(s)\big)\exp^{-(t-s)P}(-\D_v + v) \big(  \exp^{-sP}f_{0}+g(s)    \big)\ds
\end{equation*}
\begin{equation*}
\Phi_{2}(g,G)(t)= K *   \int_0^t \int \big(\Gl(s)+G(s)\big)\exp^{-(t-s)P}(-\D_v + v) \big(  \exp^{-sP}f_{0}+g(s)    \big)\dnu\ds,
\end{equation*}
and we observe that $(f,E)$ solves~\eqref{duhamel} if and only if $(g,G)=\Phi(g,G)$.  For $\s  \geq 0 $ define the norms
\begin{equation*}
\|g\|_X=\sup_{t\geq 0} \Big(\exp^{\s\kappa t}\|g(t,.)\| \Big),
\end{equation*}
\begin{equation*}
\|G\|_Y=\sup_{t\geq 0} \Big( \exp^{\s\kappa t}\|G(t,.)\|_{L^{\infty}}\Big),
\end{equation*}
define  the Banach space
$$Z:=X\times Y_{} ,\qquad \text{with}\qquad \|(g,G)\|_{Z}=\max\big(\|  g\|_X,  \|G\|_{Y_{ }}\big),$$
and   denote by $B_{\eps_0}$ the ball of $Z$ of radius $\eps_0 >0$. We will prove that if $\eps_0$ is small enough, the map $\Phi$ is a contraction of the  ball  $B_{\eps_0} \subset Z$. In the following any norm without subscript is a $L^2(\dmu\dnu)$ norm.

In a first step we prove that $\Phi$ maps $B_{\eps_0} $ into itself. Let $(g,G) \in Z$. We first estimate $\Phi_1$ and get for all $t\geq 0$
\begin{equation}
\begin{split}
\norm{\Phi_1(g,G)} \leq \int_0^t \norm{\Gl(s)+G(s)}_{L^\infty} \norm{ \exp^{-(t-s)P}(-\D_v + v) \sep{ \exp^{-sP}f_{0}+g(s) } } \ds.
\end{split}
\end{equation}
Now we use Lemma \ref{together} and we get
\begin{equation}
\begin{split}
\norm{\Phi_1(g,G)} \leq \int_0^t \exp^{-(t-s) \kappa} ((t-s)^{-1/2}  + 1) \norm{\Gl(s)+G(s)}_{L^\infty} \norm{  \exp^{-sP}f_{0}+g(s)  } \ds.
\end{split}
\end{equation}
We have on the one hand
$$
\norm{  \exp^{-sP}f_{0}+g(s)  } \leq \norm{  \exp^{-sP}f_{0}} + \norm{ g(s)  } \leq C \sep{ \eps_0 + \norm{ g  }_X}
$$
and on the other hand using Lemma \ref{estimE}
\begin{eqnarray*} %\label{380}
  \|\Gl(s)+G(s)\|_{L^{\infty}}&\leq&   \|G_l(s)\|_{L^{\infty}} + \|G(s)\|_{L^{\infty}}\nonumber \leq  C \exp^{-\s \kappa s}\|f_{0}\| + C\exp^{-\s\kappa s}\|G\|_{Y}\nonumber\\
&\leq &    C \exp^{-\s \kappa s}( \eps_0+\|G\|_{Y}).
 \end{eqnarray*}
Taking $\sigma = 1/2$ and putting the two preceding results together give
\begin{equation}
\begin{split}
\norm{\Phi_1(g,G)} & \leq C^2 \sep{ \eps_0 + \norm{ g  }_X} ( \eps_0+\|G\|_{Y}) \exp^{-\s \kappa t} \int_0^t \exp^{-\s\kappa (t-s)} ((t-s)^{-1/2}  + 1)  \ds \\
& \leq C^2 \sep{ \eps_0 + \norm{ g  }_X} ( \eps_0+\|G\|_{Y}) \exp^{-\s \kappa t},
\end{split}
\end{equation}
So that
\begin{equation*} %\label{phi1}
\begin{split}
\norm{\Phi_1(g,G)}_X & \leq C^2 \sep{ \eps_0 + \norm{ (g,G)  }_Z}^2
\end{split}
\end{equation*}
 With the same arguments and Lemma \ref{estimE} again, one easily get
\begin{equation*} %\label{phi1}
\begin{split}
\norm{\Phi_2(g,G)}_Y & \leq C^2 \sep{ \eps_0 + \norm{ (g,G)  }_Z}^2.
\end{split}
\end{equation*}
Putting these two estimates together gives
\begin{equation} \label{phi1}
\begin{split}
\norm{\Phi(g,G)}_Z & \leq C^2 \sep{ \eps_0 + \norm{ (g,G)  }_Z}^2
\end{split}
\end{equation}
Nowe let us take $\eps_0$ such that $4 C^2 \eps_0^2 \leq \eps_0$. If we suppose $\norm{ (g,G)  }_Z \leq \eps_0$ then we get that $\Phi$ maps $B_{\eps_0}$ into itself.

With exactly the  same arguments, we can also prove the contraction estimate
\begin{equation*}
\|  \Phi(g,G)-\Phi(h,H)\|_{Z} \leq C\|  (h-g,H-G)\|_{Z}\big( \eps_0 +\|(g,G)\|_{Z}+\|(h,H)\|_{Z}\big).
\end{equation*}
We do not write the details. The fixed point theorem therefore gives the existence and uniqueness
of the solution of the mollified Vlasov-Poisson-Fokker-Planck system \eqref{duhamel} with $g \in \ccc^0(\R^+, L^2)$ and $G \in \ccc^0(\R^+, L^\infty)$.
This proves point i) of  Theorem \ref{MVPFP}.

For point ii) we just have to notice that
\begin{equation*}
\|g(t)\|_{B^{}}\leq C \exp^{-\s \kappa t}\|g\|_{X} \leq \eps_0 C \exp^{-\s \kappa t}
\end{equation*}
 Similarly,
 \begin{equation*}
\|G(t)\|_{L^{\infty}}\leq C \exp^{-\s \kappa t}\|G\|_{Y}\leq \eps_0 C \exp^{-\s \kappa t}
\end{equation*}
so that writing $f(t) = \exp^{-tP} f_0 + g(t)$ and $ E(t) = K * \int f (t) \dnu$ gives the result since $\seq{f_0} = 0$.
The proof is complete. \fin
%%%%%%%%%%%%%%%%%%%%
%%%%%%%%%%%%%%%


\begin{thebibliography}{99}





\bibitem{AHL15} {Alexandre R.,  H\'erau F.,  Li W.-X.}
\it  Global hypoelliptic and symbolic estimates for the linearized Boltzmann operator
without angular cutoff,
\rm preprint

%\bibitem{AMUXY1-ARMA} AMUXY
%\bibitem{ADVW} ADVW Cancellation lemma






%\bibitem{BL04} {Bismut J.-M., Lebeau G.}
%\it le Laplacien hypoelliptique de J.-M Bismut \rm expos\'e au
%s\'eminaire X EDP octobre 2004.


%\bibitem{Bou91} {Bouchut, F.}
%\it Global weak solution of the Vlasov-Poisson system for small
%electron mass,\rm Comm. Part. Diff. \textbf{16} no. 8, 9, (1991),
%1337--1365.


\bibitem{Bou93} {Bouchut, F.}
\it Existence and uniqueness of a global smoothsolution for the
Vlasov-Poisson-Fokker-Planck system in three dimensions,
\rm J.
Funct Anal. \textbf{111}, (1993), 239--258.

%\bibitem{BD95} {Bouchut F, Dolbeault, J.}
%\it On long time asymptotics of the Vlasov-Fokker-Planck equation
%and of the Vlasov-Poisson-Fokker-Planck system with coulombic and
%Newtonian potentials, \rm Diff. Int. Eq, \textbf{8}, no 3, (1995),
%487--514.
%
%\bibitem{CSV96} { Carillo J. A., Soler J., Vasquez J.L.}
%\it Asymptotic Behavior and selfsimilarity for the three
%dimentional Vlasov-Poisson-Fokker-Planck System, \rm Journ. Funct.
%Ana. 141, (1996), 99--132.

%\bibitem{CTW16} Carraptoso Tristani Wu Landau



\bibitem{CSV92} { Coulhon T., Saloff-Coste L. Varopoulos N. Th.}
\it Analysis and geometry on groups, \rm Cambridge Tracts in
Math., Cambridge University press (1992).




\bibitem{Deg86} { Degond, P.}
\it Global existence of smooth solutions for the
Vlasov-Poisson-Fokker-Planck equation in 1 and 2 space dimensions
\rm Ann. scient. Ec. Norm. Sup., 4ème s\'erie,  19, (1986),
519--542.



%\bibitem{DV01} {Desvillettes L. and Villani C.}
%\it On the trend to global equilibrium in spatially inhomogeneous
%systems. Part I: the linear Fokker-Planck equation. \rm Comm. Pure
%Appl. Math. 54, 1 (2001), 1-42


\bibitem{DHMS16} {Dolbeault J., H\'erau F., Mouhot C., Schmeiser C.}
\it Hypocoercivity for linear kinetic equations with several conservation laws,
\rm in progress

\bibitem{DMS15} {Dolbeault J., Mouhot C., Schmeiser C.}
\it Hypocoercivity for linear kinetic equations conserving mass
\rm Trans. Am. Math. soc. 367, No 6, 2015, 3807--3828


%\bibitem{DV04} {Desvillettes L. and Villani C.}
%\it On the trend to global equilibrium for spatially inhomogeneous
%systems: the Boltzmann equation. \rm Invent. Math. 105, 1 (2004).


%\bibitem{Dol91}  {Dolbeault, J.}
%\it Stationary states in plasma physics: Maxwellian solutions of
%the Vlasov-Poisson system. \rm Math. Mod. Meth. Appli. Sci.
%\textbf{1} no. 2,(1991), 183--208.



%\bibitem{Dol99}  {Dolbeault, J.}
%\it Free energy and solutions of the Vlasov-Poisson-Fokker-Planck
%system: external potential and confinment (large time behavior and
%steady states). \rm J. Math.pures. Appl., \textbf{78},(1999),
%121--157.


%\bibitem{Dre87}  {Dressler, K.}
%\it Stationary solutions of the Vlasov-Fokker-Planck equations.
%\rm Math. Meth.Appl. Sci., \textbf{9},(1987), 169--176.


%\bibitem{Dre90}  {Dressler, K.}
%\it Steady states in plasma physics -- The Vlasov-Fokker-Planck
%equation \rm Math. Meth.Appl. Sci., \textbf{12} no. 6,(1990),
%471--487.


\bibitem{DHL16} {Dujardin G., H\'erau F., Laffite P.},
\it Coercivity and hypocoercivity for the discrete Fokker-Planck equation.
\rm in Progress

%\bibitem{EH00} {Eckmann J.-P., Hairer M.}
%\it{Non-equilibrium statistical mechanics of strongly anharmonic
%chains of oscillators}, \rm Comm. Math. Phys. 212 , no. 1,
%105--164 (2000).





%\bibitem{EH03} {Eckmann J.-P., Hairer M.}
%\it{spectral properties of hypoelliptic operators}, \rm Commun.
%Math. Phys. \textbf{235} no.2, (2003), 233--253.
%
%\bibitem{EPR99} {Eckmann J.-P., Pillet, C.-A. et Rey-Bellet L.}
%\it{Non-equilibrium statistical mechanics of anharmonic chains
%coupled to two heat bath at different temperature}, \rm {Comm.
%Math. Phys.}, \textbf{201}, no.3, (1999),  657--697.


\bibitem{FS86} {Fefferman C., Sanchez-Calle S.}
\it{Fundamental solutions for second order differential
operators}, \rm {Ann. Math.}, \textbf{124},  (1986), 247-272.



%\bibitem{Guo02} {Guo, yan}
%\it{The Landau equation in a periodic box}, \rm {Commun. Math.
%Phys.}, \textbf{231}, (2002), 391-434.


%\bibitem{GSZ96} Glassey R., Schaeffer J., Zheng Y.
%\it  Steady states of the Vlasov-Poisson-Fokker-Planck System \rm
%J.  Math. Ana. Appl. \textbf{202}, (1996), 1058--1075.


%\bibitem{GL89} Gogny D. Lions P.L.
%\it  Sur les \'etats d'\'equilibre pour les densit\'es \'el\'ectroniques
%dans les plasma \rm RAIRO mod\'el.,  Math. Ana. Num. \textbf{23} no.
%1, (1989), 137--153.



\bibitem{HelN04} Helffer B., Nier F.
\newblock  Hypoelliptic estimates and spectral theory for Fokker-Planck operators and Witten Laplacians.
\newblock {\em Lecture Notes in Mathematics}, 1862. Springer-Verlag, Berlin, 2005. x+209 pp.


\bibitem{Her06}
 {H\'erau, F.},
 \newblock  \it Hypocoercivity and exponential time decay for the linear inhomogeneous relaxation Boltzmann equation,
 \newblock \rm {Asymptot. Anal. 46} (2006), 349--359.


 \bibitem{Her07}
H\'erau F.,
\newblock  \it Short and long time behavior of the Fokker-Planck equation in a confining potential and applications.
\newblock \rm J. Funct. Anal.,  \textbf{244}, no. 1, 95--118 (2007).


\bibitem{HN04}
H{\'e}rau F. and Nier F.
\newblock{\em  Isotropic hypoellipticity and trend to
equilibrium for the
  Fokker-Planck equation with high degree potential}.
\newblock { Arch. Ration. Mech. Anal.}, 171(2),
 151--218, 2004. announced in
Actes colloque EDP Forges-les-eaux, 12p.,  (2002).


%\bibitem{HSS05} {H\'erau F., Sj\"ostrand J., Stolk, C. },
%\it Semiclassical analysis for the Kramers-Fokker-Planck equation,
%\rm to appear in Comm. Part. Diff. Eq.,  (2005)

\bibitem{HT14} {H\'erau F., Thomann L.},
\it On global existence and trend to the equilibrium for the Vlasov-Poisson-Fokker-Planck
system with exterior confining potential
\rm J. F. A.  271 (5), 1301--1340

\bibitem{HTT16a} {H\'erau F.,  Tonon D.,  Tristani, I.},
\it Short time diffusion properties of inhomogeneous kinetic equations with fractional collision kernel, \rm in progress

\bibitem{HTT16b} {H\'erau F.,  Tonon D.,  Tristani, I.},
\it Cauchy theory and exponential stability for inhomogeneous Boltzmann equation for hard potentials without cut-off, \rm in progress



\bibitem{Hor67} {H{\"o}rmander L.},
\it  {Hypoelliptic second order differential equations}, \rm Acta
Math. \textbf{119}, (1967), 147--171.


\bibitem{Hor95}
H\"ormander L.,
\it Symplectic classification of quadratic forms, and general Mehler formulas. \rm \newblock{ Mathematische Zeitschrift},  \textbf{219} (3), 413--450, 1995



%\bibitem{Hor85} {H{\"o}rmander, L.},
%\it  {The analysis of linear partial differential operators
%{I},{I}{I}{I}}, \rm  {Springer-Verlag},  {Berlin}, (1985).


%\bibitem{Kag01} {Kagei Y.}
%\newblock {\em Invariant manifolds and long-time
%asymptotics for the Vlasov-Poisson-Fokker-Planck equation}. \rm
%SIAM J. Math. Anal. \textbf{33} no. 2, (2001), 489--507.


\bibitem{Koh73} {Kohn J. J.}
\newblock {\em Pseudodifferential operators and hypoellipticity}. \rm
Proc. Symp. Pure Math. 23, AMS    61-69, 1973.


\bibitem{GMM13} {Gualdani M.-P., Mouhot C., Mischler S.} \it Factorization for non-symmetric operators and exponential H-theorem. \rm preprint

 \bibitem{MN06}
 { Mouhot, C., Neumann, L.}
 \newblock \it Quantitative perturbative study of convergence to equilibrium for collisional kinetic models in the torus.
 \newblock \rm { Nonlinearity 19} (2006), 969--998.

%\bibitem{Lio61} {Lions J. L.}
%\newblock {\em Equations diff\'erentielles op\'erationelles
%et problèmes aux limites}.
%\newblock Springer, Berlin, 1961.




%\bibitem{Nel67}
%Nelson E.
%\newblock {\em Dynamical theories of {B}rownian motion}.
%\newblock Princeton University Press, Princeton, N.J., 1967.




%\bibitem{OV90} {O'Dwyer B.P. and Victory H.P.}
%\newblock {\em On classical solutions of
%Vlasov-Poisson-Fokker-Planck systems}.
%\rm Indiana Univ.  Math. J. \textbf{39}, (1990), 105--157.



%\bibitem{OS00} {Ono K. and Strauss W.}
%\it Regular solutions of the Vlasov-Poisson-Fokker-Planck system.
%\rm Dis. Cont. Dyn. Syst. \textbf{6} no. 4, (2000), 751--772.







%\bibitem{Ono01} {Ono K.},
%\it Asymptotic behavior for the Vlasov-Poisson-Fokker-Planck
%System and the Collision-less Vlasov-Poisson system. \rm Non. Ana.
%\textbf{47}, (2001), 2539--2550.



%\bibitem{Paz83}
%Pazy A.
%\newblock{\em Semigroups of linear operators and
%applications to partial differential equations}.
%\newblock Springer-Verlag, Berlin, second edition, 1983.
%\newblock Applied Mathematical Sciences.




%\bibitem{Paz83}
%Pazy A.
%\newblock{\em Semigroups of linear operators and
%applications to partial differential equations}.
%\newblock Springer-Verlag, Berlin, second edition, 1983.
%\newblock Applied Mathematical Sciences.

\bibitem{PZ16} Porretta A., Zuazua E.
\it Numerical hypocoercivity for the kolmogorov equation,
\rm to appear in Math. of Comp., 2016


%\bibitem{RW92}{ G. Rein and J. Weckler}
%\newblock{\em Generic global classical solutions of the
%Vlasov-Fokker-Planck-Poisson system in three dimensions}. \rm J.
%Diff. Eq. \textbf{99}, (1992), 59--77.


%\bibitem{Ris89}
%Risken H.
%\newblock {\em The {F}okker-{P}lanck equation}.
%\newblock Springer-Verlag, Berlin, second edition, 1989.
%\newblock Methods of solution and applications.


%\bibitem{Sol97}
%Soler J.
%\newblock {\em Asymptotic behavior for the
%Vlasov-Poisson-Fokker-Planck system}.
%\newblock Non. lin. Ana. th. Meth. Appl., \textbf{30} no. 8,
% proc. 2nd World Congress of Nonlin. Ana, (1997), 5217--5228.


%\bibitem{Ste70}
%Stein E.M.
%\newblock {\em Singular integrals and differentiability of
%functions}.
%\newblock Princeton University Press, Princeton, 1970.



%\bibitem{Tal99}
%Talay D.
%\newblock \it Approximation of invariant
% measures of nonlinear {H}amiltonian and
%  dissipative stochastic differential equations. \rm
%\newblock In C.~Soize R.~Bouc, editor,
%{\em Progress in Stochastic Structural
%  Dynamics}, volume 152 of {\em Publication du L.M.A.-C.N.R.S.},
%   (1999), 139--169.



%\bibitem{Vil02}
% Villani C.
%\it A review of mathematical topics in collisional kinetic theory.
%\rm Handbook of Fluid Mechanics S. Friedlander and D. Serre Eds.
%(2002).

\bibitem{Tri15}  Tristani I.,  \it Fractional Fokker-Planck equation. \rm Commun. Math. Sci. 13, 5 (2015), 1243--1260

 \bibitem{Vil09}
Villani C.,
\newblock \it Hypocoercivity.
\newblock \rm { Mem. Amer. Math. Soc.} 202 (2009), no. 950, iv+141 pp.



%\bibitem{Won99}
% Wong M.W.
%\it An introduction to pseudo-differential operators. \rm World
%Scientific, second edition, 1999.



% \bibitem{villani}
%C. Villani.
%\newblock Hypocoercivity.
%\newblock {\em Mem. Amer. Math. Soc.} 202 (2009), no. 950, iv+141 pp.





 \end{thebibliography}
\end{document}